\documentclass[reqno]{gokova}
\usepackage{epsfig,graphicx,amscd}

\begin{document}
\def\E{\ifmmode{\mathbb E}\else{$\mathbb E$}\fi} 
\def\N{\ifmmode{\mathbb N}\else{$\mathbb N$}\fi} 
\def\R{\ifmmode{\mathbb R}\else{$\mathbb R$}\fi} 
\def\Q{\ifmmode{\mathbb Q}\else{$\mathbb Q$}\fi} 
\def\C{\ifmmode{\mathbb C}\else{$\mathbb C$}\fi} 
\def\H{\ifmmode{\mathbb H}\else{$\mathbb H$}\fi} 
\def\Z{\ifmmode{\mathbb Z}\else{$\mathbb Z$}\fi} 
\def\P{\ifmmode{\mathbb P}\else{$\mathbb P$}\fi} 
\def\T{\ifmmode{\mathbb T}\else{$\mathbb T$}\fi} 
\def\SS{\ifmmode{\mathbb S}\else{$\mathbb S$}\fi} 
\def\DD{\ifmmode{\mathbb D}\else{$\mathbb D$}\fi} 

\renewcommand{\a}{\alpha}
\renewcommand{\b}{\beta}
\renewcommand{\d}{\delta}
\newcommand{\D}{\Delta}
\newcommand{\e}{\varepsilon}
\newcommand{\g}{\gamma}
\newcommand{\G}{\Gamma}
\newcommand{\la}{\lambda}
\newcommand{\La}{\Lambda}
\newcommand{\n}{\nabla}
\newcommand{\var}{\varphi}
\newcommand{\s}{\sigma}
\newcommand{\Sig}{\Sigma}
\renewcommand{\t}{\tau}
\renewcommand{\th}{\theta}
\renewcommand{\O}{\Omega}
\renewcommand{\o}{\omega}
\newcommand{\z}{\zeta}

\newcommand{\ben}{\begin{enumerate}}
\newcommand{\een}{\end{enumerate}}
\newcommand{\be}{\begin{equation}}
\newcommand{\ee}{\end{equation}}
\newcommand{\bea}{\begin{eqnarray}}
\newcommand{\eea}{\end{eqnarray}}
\newcommand{\bc}{\begin{center}}
\newcommand{\ec}{\end{center}}

\newtheorem{thm}{Theorem}[section]
\newtheorem{cor}[thm]{Corollary}
\newtheorem{lem}[thm]{Lemma}
\newtheorem{prop}[thm]{Proposition}
\newtheorem{ax}{Axiom}
\newtheorem{conj}[thm]{Conjecture}

\theoremstyle{definition}
\newtheorem{defn}{Definition}[section]

\theoremstyle{remark}
\newtheorem{rem}{\rm\bfseries{Remark}}[section]
\newtheorem*{notation}{Notation}

\newtheorem{ques}{\rm\bfseries{Question}}[section]
\newtheorem{cons}[rem]{\rm\bfseries{Construction}}
\newtheorem{exm}[rem]{\rm\bfseries{Example}}






\def\C{\mathbb C}
\def\H{\mathbb H}
\def\I{\mathbb I}
\def\P{\mathbb P}
\def\Q{\mathbb Q}
\def\R{\mathbb R}
\def\T{\mathbb T}
\def\Z{\mathbb Z}
\def\D{\mathbb D}

\def\Ga{\alpha}
\def\Gb{\beta}
\def\Ge{\varepsilon}
\def\Gg{\gamma}
\def\GG{\Gamma}
\def\Gd{\delta}
\def\GD{\Delta}
\def\GH{\aleph}
\def\GL{\Lambda}
\def\Gl{\lambda}
\def\Go{\omega}
\def\GO{\Omega}
\def\Gr{\varrho}
\def\GP{\Phi}
\def\Gp{\phi}
\def\Gvp{\varphi}
\def\GS{\Sigma}
\def\Gs{\sigma}
\def\Gvs{\varsigma}
\def\Gt{\tau}

\def\cA{\mathcal A}
\def\cB{\mathcal B}
\def\cC{\mathcal C}
\def\cD{\mathcal D}
\def\cE{\mathcal E}
\def\cF{\mathcal F}
\def\cG{\mathcal G}
\def\cI{\mathcal I}
\def\cJ{\mathcal J}
\def\cL{\mathcal L}
\def\cM{\mathcal M}
\def\cN{\mathcal N}
\def\cO{\mathcal O}
\def\cP{\mathcal P}
\def\cR{\mathcal R}
\def\cQ{\mathcal Q}
\def\cS{\mathcal S}
\def\cT{\mathcal T}
\def\cU{\mathcal U}

\def\fA{\mathfrak A}
\def\fB{\mathfrak B}
\def\fc{\mathfrak c}
\def\fd{\mathfrak d}
\def\fe{\mathfrak e}
\def\fF{\mathfrak F}
\def\fM{\mathfrak M}
\def\fr{\mathfrak r}
\def\fU{\mathfrak U}

\def\dd{\partial}
\def\smin{\setminus}
\def\fix{\mathrm{Fix}}
\def\Conj{\mathrm{Perm}}
\def\conj{conj}
\def\emp{\emptyset}

\def\Nn{1,\cdots,n}
\def\im{\sqrt{-1}}
\def\Aut{\operatorname{Aut}}
\def\one{1\!\! 1}
\def\rev{\curvearrowright}
\def\projc{\C\P^{1}}
\def\projr{\R\P^{1}}




\def\epsffile#1{\psfig{file=#1,silent=}}
\def\eq#1{\underset{#1}{=}}
\def\fig#1{\vcenter{\psfig{figure=#1,silent=}}}
\def\centerfig#1{\centerline{\psfig{figure=#1,silent=}}}


\newcommand{\orT}[1]{(#1)}
\def\tree{\mathcal Tree}
\def\atree{\mathcal ATree(\Gs)}
\def\ptree{\mathcal PTree(\Gs)}
\def\otree{\mathcal UTree(\Gs) }

\newcommand{\ve}[1]{\mathbf{#1}}

\newcommand{\real}[1]{\R #1}

\def\curve{(\GS;{\mathbf p})}
\def\curvef{(\overline{\GS}; \Gs({\mathbf p}))}

\newcommand{\cmod}[1]{\overline{M}_{#1}}
\newcommand{\cmodo}[1]{M_{#1}}
\newcommand{\umod}[1]{\overline{\cU}_{#1} }
\newcommand{\tumod}[1]{\widetilde{\cU}_{#1} }
\newcommand{\umodo}[1]{\cU_{#1} }

\newcommand{\rmod}[1]{\R \overline{M}_{(#1)} }
\newcommand{\rmodo}[1]{\R M_{(#1)} }

\newcommand{\cover}[1]{\R \widetilde{M}_{(#1)}}
\newcommand{\ocover}[1]{\R \widetilde{M}^{o}_{#1}}

\newcommand{\cspp}[1]{{\widetilde{Conf}}_{(#1)}}
\newcommand{\csp}[1]{Conf_{(#1)}}

\newcommand{\cspq}[1]{\widetilde{C}_{(#1)}}
\newcommand{\csq}[1]{C_{(#1)}}
\newcommand{\csqc}[1]{\overline{C}_{(#1)}}

\newcommand{\dsq}[1]{B_{(#1)}}
\newcommand{\dsqc}[1]{\overline{B}_{(#1)}}

\newcommand{\form}[1]{\Go_{(#1)}}
\newcommand{\ori}[1]{ [ \Go_{(#1)} ] }
\newcommand{\orie}[1]{[\GO_{(#1)}]}

\newcommand{\konj}[1]{\overline{#1}}

\newcommand{\ove}[1]{\widetilde{#1}}

\newcommand{\call}[1]{\left[ \overline{D}_{#1} \right]}

\newcommand{\open}[1]{U_{(#1)} }
\newcommand{\res}[2]{\mathfrak{res}_{#1}^{#2}}

\setcounter{page}{1}
\volume{13}

\title[Moduli of pointed real curves of genus zero]{On
moduli of pointed real curves of genus zero}
\author[CEYHAN]{\"Ozg\"ur Ceyhan}

\thanks{}

\address{Centre de Recherches Math\'ematiques, Universit\'e de Montr\'eal, Montr\'eal, Canada}
\email{ceyhan@crm.umontreal.ca}

\begin{abstract}
We introduce the moduli space $\rmod{2k,l}$ of pointed real curves
of genus zero and give its natural stratification. The strata of
$\rmod{2k,l}$ correspond to real curves of genus zero with different
degeneration types and are encoded by trees with certain
decorations. By using this stratification, we calculate the first
Stiefel-Whitney class of $\rmod{2k,l}$ and construct the orientation
double cover $\cover{2k,l}$ of $\rmod{2k,l}$.
\end{abstract}

\maketitle

\section{Introduction}
The moduli space $\cmod{n}$ of stable $n$-pointed (complex) curves
of genus zero has been extensively studied as one of the fundamental
models of moduli problems in algebraic geometry (see
\cite{ke,knu,km1,km2,m}). The moduli space $\cmod{n}$ carries a set
of anti-holomorphic involutions, whose fixed point sets are the
moduli spaces of pointed real curves of genus zero. These moduli
spaces parameterize the isomorphism classes of pointed curves of
genus zero with a real structure. For each of these spaces, a
certain set of labeled points stays in the real parts of the curves
while other pairs of labeled points are conjugated by the real
structures of the curves.

The moduli spaces of pointed real curves have recently attracted
attention in various contexts such as multiple $\zeta$-motives
\cite{gm}, representations of quantum groups \cite{ehkr,hk} and
Welschinger invariants \cite{w1,w2}.

The aim of this work is to explore the topological properties of the
moduli spaces of pointed real curves of genus zero. Hence, we first
introduce a natural combinatorial stratification of the moduli
spaces of pointed real curves of genus zero through the
stratification of $\cmod{n}$. Each stratum is determined by the
degeneration type of the real curve. They are identified with the
product of the spaces of real point configurations in the projective
line $\projc$ and the moduli spaces $\cmod{m}$. The degeneration
types of the pointed real curves are encoded by trees with
corresponding decorations. Secondly, we calculate the first
Stiefel-Whitney classes of the moduli spaces in terms of their
stratifications. The moduli spaces of pointed real curves are not
orientable for $n \geq 5$ and the set of labeled real points of
real curves is not
empty. We construct the orientation double covers of the moduli
spaces for the non-orientable cases. The double covering in this
work significantly differs from the `double covering' in the recent
literature on open Gromov-Witten invariants and moduli spaces of
pseudoholomorphic discs (see \cite{fu,liu}): Our double covering has
no boundaries which suits better for the use of intersection theory.

This paper is organized as follows. Section \ref{ch_n_curves} contains 
a brief overview of facts about the moduli space $\cmod{n}$. In Section
\ref{ch_real_curves}, we introduce real structures on $\cmod{n}$ and
real parts $\rmod{2k,l}$ as moduli spaces of $(2k,l)$-pointed real
curves. The following section, the stratification of $\rmod{2k,l}$
is given according to the degeneration types of pointed real curves
of genus zero. In Section \ref{sec_steifel}, the first
Stiefel-Whitney class of $\rmod{2k,l}$ is calculated by using the
stratification given in Section \ref{sec_stratification_sp}. Then in
Section \ref{sec_coverofmoduli}, we construct the orientation double
coverings $\cover{2k,l} \to \rmod{2k,l}$.

In this paper, the genus of the curves is zero except when the
contrary is stated explicitly. Therefore, we omit mentioning the
genus of the curves.

\section{Pointed complex curves and their moduli}
\label{ch_n_curves}

This section reviews the basic facts on pointed complex curves of
genus zero and their moduli space.

\subsection{Pointed curves and their trees}
\label{sec_n_curves}

\begin{defn}
\label{def_n-curve}%
An {\it $n$-pointed curve} $\curve$ is a connected complex algebraic
curve $\GS$ with   distinct, smooth, {\it labeled points}
$\mathbf{p} = (p_{1},\cdots,p_{n}) \subset \GS$, satisfying the
following conditions:
\begin{itemize}
\item $\GS$ has only nodal singularities.%
\item The arithmetic genus of $\GS$ is equal to zero.
\end{itemize}
A {\it family of $n$-pointed curves} over a complex  manifold $S$ is
a proper, holomorphic map $\pi_S:\cU_S \to S$ with $n$ sections
$p_1,\cdots,p_n$ such that each geometric fiber
$(\GS(s);\mathbf{p}(s))$ is an $n$-pointed curve.

Two such  curves,  $\curve$ and $(\GS';\ve{p}')$, are {\it isomorphic}
if there exists a bi-holomorphic equivalence $\Phi: \GS \to \GS'$
mapping $p_{i}$ to $p'_{i}$.

An $n$-pointed curve is {\it stable} if its automorphism group is
trivial (i.e., on each irreducible component, the number of singular
points plus the number of labeled points is at least three).
\end{defn}

\subsubsection{Graphs}

\begin{defn}
\label{def_graph}%
A {\it graph} $\GG$ is a collection of finite sets of {\it vertices}
$V_{\GG}$ and {\it flags (or half edges)} $F_{\GG}$ with a boundary
map $\dd_{\GG}: F_{\GG} \to V_{\GG}$ and an involution $j_{\GG}:
F_{\GG} \to F_{\GG}$ ($j_{\GG}^{2} =id$). We call $E_{\GG} =
\{(f_{1},f_{2}) \in F_{\GG}^{2} \mid f_{1}= j_{\GG} f_{2}\ \& \
f_{1} \not= f_{2} \}$ the set of {\it edges}, and $T_{\GG} = \{f \in
F_{\GG} \mid f = j_{\GG} f\}$ the set of {\it tails}. For a vertex
$v \in V_{\GG}$, let $F_{\GG}(v) =\dd^{-1}_{\GG}(v)$ and $|v| =
|F_{\GG}(v)|$ be the {\it valency} of $v$.
\end{defn}

We think of a graph $\GG$ in terms of  its following {\it geometric
realization} $||\GG||$: Consider the disjoint union of closed
intervals $\bigsqcup_{f_i \in F_{\GG}} [0,1]  \times f_i$ and
identify $(0, f_i)$ with $(0, f_j)$ if $\dd_\GG f_i = \dd_\GG f_j$,
and identify $(t, f_i)$ with $(1-t, j_\GG f_i)$ for $t \in ]0,1[$
and $f_i \ne j_\GG f_i$. The geometric realization of $\GG$ has a
piecewise linear structure.

\begin{defn}
\label{def_tree}%
A {\it tree} $\Gg$ is a graph whose geometric realization is
connected and simply-connected. If $|v|>2$ for all $v \in V_{\Gg}$,
then such a tree is called {\it stable}.
\end{defn}

We associate a subtree $\Gg_{v}$ for each vertex $v \in V_{\Gg}$
which is given by $V_{\Gg_{v}} = \{v\}, F_{\Gg_{v}} = F_{\Gg}(v)$,
$j_{\Gg_{v}} =id$, and  $\dd_{\Gg_{v}} = \dd_\Gg$.

\begin{defn}
\label{def_tree_morphism}%
Let $\Gg$ and $\Gt$ be  trees with $n$ tails.  A {\it morphism}
between these trees $\Gp: \Gg\to \Gt$ is a pair of maps $\Gp_{F}:
F_{\Gt} \to F_{\Gg}$ and $\Gp_{V}: V_{\Gg} \to V_{\Gt}$  satisfying
the following conditions:
\begin{itemize}
\item $\Gp_{F}$ is injective and $\Gp_{V}$ is surjective. %
\item The following diagram commutes
$$
\begin{CD}
F_{\Gg}  @> {\dd_{\Gg}} >>   V_{\Gg} \\
@A{\operatorname{\Gp}_F}AA      @VV{\operatorname{\Gp}_V}V \\
F_{\Gt}   @> {\dd_{\Gt}}>> V_{\Gt}.
\end{CD}
$$
\item $\Gp_{F} \circ j_{\Gt} = j_{\Gg} \circ \Gp_{F}$.
\item $\Gp_{T}:={\Gp_{F} |}_{T}$ is a bijection.
\end{itemize}
An {\it isomorphism} $\Gp: \Gg \to \Gt$ is a morphism where
$\Gp_{F}$ and $\Gp_{V}$ are bijections. We denote the isomorphic
trees by $\Gg \approx \Gt$.
\end{defn}

Each morphism induces a piecewise linear map on geometric
realizations.

\begin{lem}
\label{lem_tree_iso}%
Let $\Gg$ and $\Gt$ be stable trees with $n$ tails. Any isomorphism
$\phi: \Gg \to \Gt$ is uniquely defined by its restriction on tails
$\phi_{T}: T_{\Gt} \to T_{\Gg}$.
\end{lem}

\begin{proof} Let $\phi, \varphi: \Gg \to \Gt$ be two isomorphisms
such that their restriction on tails are the same. Consider the path
$_{f_1}P_{f_2}$ in  $||\Gg||$ that connects a pair of tails
$f_1,f_2$. The automorphism $\varphi^{-1} \circ \phi$ of $\Gg$ maps
$_{f_1}P_{f_2}$ to itself; otherwise, the union of the
$_{f_1}P_{f_2}$ and its image $\varphi^{-1} \circ
\phi(_{f_1}P_{f_2})$ gives a loop in $||\Gg||$, which contradicts
simply-connectedness. Moreover, the restriction of $\varphi^{-1}
\circ \phi$ to the path $_{f_1}P_{f_2}$ is the identity map since it
preserves distances of vertices to tails $f_1,f_2$. This follows
from the compatibility of the automorphism $\varphi^{-1} \circ \phi$
with $\dd_\Gg$ and $j_\Gg$.

The geometric realization $||\Gg||$ of $\Gg$ can be covered by paths
that connects pairs of tails of $\Gg$. We conclude that the
automorphism $\varphi^{-1} \circ \phi$ is the identity since it is
the identity on every such path.
\end{proof}

There are only finitely many isomorphisms classes of stable trees
whose set of tails is $T_{\Gg} = \{1,\cdots,n\}$. We call the
isomorphism classes of such trees {\it $n$-trees}. We denote the set
of all $n$-trees by $\tree$.

\subsubsection{Dual trees of pointed curves}
Let $\curve$ be an $n$-pointed curve and $\eta: \hat{\GS} \to \GS$
be its normalization. Let $(\hat{\GS}_v;\hat{\mathbf{p}}_v)$ be the
following $|v|$-pointed stable curve: $\hat{\GS}_{v}$ is a component
of $\hat{\GS}$, and $\hat{\mathbf{p}}_v$ is the set of points
consisting of the preimages of {\it special} (i.e.  labeled and
nodal) points on $\GS_{v} := \eta (\hat{\GS}_v)$. The points
$\hat{\mathbf{p}}_v = (p_{f_{1}},\cdots, p_{f_{|v|}})$ on
$\hat{\GS}_v$ are ordered by the flags $f_* \in F_{\Gt}(v)$.

\begin{defn}
\label{def_dual_tree}%
The {\it dual tree} $\Gg$ of an  $n$-pointed curve $\curve$ is the
tree consisting of following data:
\begin{itemize}
\item $V_{\Gg}$ is the set of  components of $\hat{\GS}$.
\item $F_{\Gg}$ is the set consisting of the preimages of special points.
\item $\dd_{\Gg}: f \mapsto v$ if and only if $p_f \in \hat{\GS}_v$.
\item $j_{\Gg}: f \mapsto f$ if and only if $p_f$ is a labeled point, and
$j_{\Gg}: f_1 \mapsto f_2$ if and only if $p_{f_1} \in
\hat{\GS}_{v_1}$ and $p_{f_2} \in  \hat{\GS}_{v_2}$ are the
preimages of  the nodal point $\GS_{v_1} \cap \GS_{v_2}$.
\end{itemize}
\end{defn}

\begin{lem} \label{lem_bihol}
Let $\Phi$ be an isomorphism between the $n$-pointed stable curves
$\curve$ and $ (\GS';\mathbf{p}')$.

(i) $\Phi$ induces an isomorphism $\phi$ between their dual trees
$\Gg,\Gt$.

(ii) $\Phi$  is uniquely defined by its restriction on labeled
points.
\end{lem}

\begin{proof} (i)  The result follows from the decomposition of $\GP$ into
its restriction to each irreducible component and the Def.
\ref{def_tree_morphism}.

(ii) Due to Lemma \ref{lem_tree_iso}, the isomorphism $\Gp: \Gg \to
\Gt$ is determined by the restriction of $\GP$ to the labeled
points. The isomorphism $\Gp$ determines which component of $\GS$ is
mapped to which component of $\GS'$ as well as the restriction of
$\GP$ to the special points. Each component of $\GS$ is rational and
has at least three special points. Therefore, the restriction of
$\GP$ to a component is uniquely determined by the images of the
three special points.
\end{proof}

\subsection{Deformations of pointed curves}
\label{sec_deform}%

Let $\Gg$  be the dual tree of $\curve$ and $\hat{\GS} \to \GS$ be
the normalization. Let $(\hat{\GS}_v;\hat{\mathbf{p}}_v)$ be the
following $|v|$-pointed stable curve corresponding to the
irreducible component $\GS_v$ of $\curve$. Let $\GO^1_{\GS }$ be the
sheaf of K\"ahler differentials.

The infinitesimal  deformations of a nodal curve $\GS$ with  divisor
$D_{\mathbf{p}}= p_1+\cdots+p_n$ is canonically identified with the
complex vector space
\begin{eqnarray} \label{eqn_def}
Ext^{1}_{\cO_{\GS}} (\GO^1_{\GS}(D_{\mathbf{p}}),\cO_{\GS}),
\end{eqnarray}
and the obstruction  lies in
\begin{eqnarray}
Ext^{2}_{\cO_{\GS}} (\GO^1_{\GS}(D_{\mathbf{p}}),\cO_{\GS}).
\nonumber
\end{eqnarray}
In this case, it is known that there are no obstructions (see, for
example \cite{liu} or \cite{hm}).

The space of infinitesimal  deformations is the tangent space of the
space of deformations at $\curve$. It can be written explicitly in
the following form:
\begin{eqnarray} \label{eqn_deform}
\bigoplus_{v \in V_{\Gg}} H^{1} (\hat{\GS}_{v},T_{\hat{\GS}_v} (-
D_{\hat{\ve{p}}_v} )) \oplus \bigoplus_{ (f_e,f^e) \in E_{\Gg}}
T_{p_{f_e}} \hat{\GS} \otimes T_{p_{f^e}}\hat{\GS}.
\end{eqnarray}
The first part corresponds to the equisingular deformations of $\GS$
with the divisor $D_{\hat{\mathbf{p}}_v} = \sum_{f_i \in F_{\Gg}(v)}
p_{f_i}$, and the second part corresponds to the smoothing of nodal
points $p_e$ of the edges $e=(f_e,f^e)$ (see \cite{hm}).

\subsubsection{Combinatorics of degenerations}
\label{sec_degen}%

Let $\curve$ be an $n$-pointed curve with the dual tree $\Gg$.
Consider the deformation of a nodal point of  $\curve$. Such a
deformation of $\curve$ gives a {\it contraction} of an edge of
$\Gg$: Let $e=(f_{e},f^{e}) \in E_{\Gg}$ be the edge corresponding
to the nodal point and $\dd_{\Gg}(e) = \{v_{e},v^{e}\}$, and
consider the equivalence relation $\sim$ on the set of vertices,
defined by: $v \sim v$ for all $v \in V_{\Gg} \smin
\{v_{e},v^{e}\}$, and  $v_{e} \sim v^{e}$. Then, there is an
$n$-tree $\Gg/e$ whose vertices are $V_{\Gg}/\sim$ and whose flags
are $F_{\Gg} \smin \{f_e,f^e\}$.  The boundary map and involution of
$\Gg/e$ are the restrictions of $\dd_{\Gg}$ and $j_{\Gg}$.

We use the notation $\Gg < \Gt$ to indicate that $\Gt$ is obtained
by contracting some edges of $\Gg$.

\subsection{Stratification of the moduli space $\cmod{n}$}
\label{stratacomplex}%

The moduli space $\cmod{n}$ is the space of  isomorphism classes of
$n$-pointed stable curves. This space is stratified according to
degeneration types of $n$-pointed stable curves which are given by
$n$-trees. The principal stratum $\cmodo{n}$ corresponds to the
one-vertex $n$-tree and is the quotient of the product
$(\projc)^{n}$ minus the diagonals $\GD = \bigcup_{k<l} \{(p_{1},
\cdots,p_{n})| p_{k} = p_{l} \}$ by $Aut(\projc) = PSL_2(\C)$.

\begin{thm}[Knudsen \& Keel, \cite{knu,ke}]

(i) For any $n \geq 3$,  $\cmod{n}$ is a smooth projective algebraic
variety of (real) dimension $2n-6$.

(ii)  Any family of $n$-pointed stable curves over $S$ is induced by
a unique morphism $\kappa: S \to \cmod{n}$. The universal family of
curves $\umod{n}$ of $\cmod{n}$ is isomorphic to $\cmod{n+1}$.

(iii) For any $n$-tree $\Gg$, there exists a quasi-projective
subvariety
 $D_{\Gg} \subset \cmod{n}$ parameterizing the curves whose dual tree
is given by $\Gg$. $D_{\Gg}$ is isomorphic to $\prod_{v \in V_{\Gg}}
\cmodo{|v|}$. Its (real) codimension is $2|E_{\Gg}|$.

(iv) $\cmod{n}$ is stratified by pairwise disjoint subvarieties
$D_{\Gg}$. The closure of any stratum $D_{\Gg}$  is stratified by
$\{D_{\Gg'} \mid \Gg' \leq \Gg \}$.
\end{thm}

\subsubsection{Examples}
\label{exa_mod}%

(i) For $n < 3$, $\cmod{n}$ is empty due to the definition of
$n$-pointed stable curves. $\cmod{3}$ is simply a point, and its
universal curve $\umod{3}$ is $\projc$ endowed with three points.

\begin{figure}[htb]
\centerfig{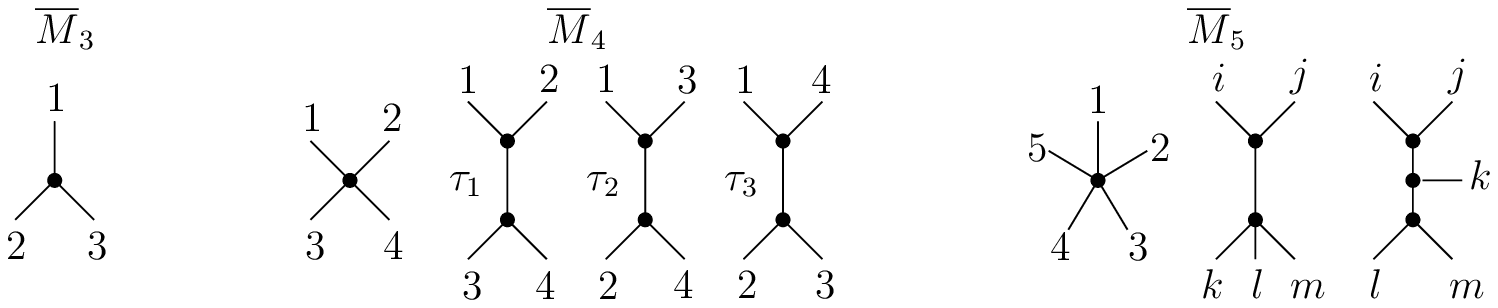,height=1.in} %
\caption{Dual trees encoding the strata of $\cmod{n}$ for $n=3,4$, and $5$.}%
\label{moduli}
\end{figure}

(ii) The moduli space $\cmod{4}$ is $\projc$ with three points.
These points $D_{\Gt_{1}},D_{\Gt_{2}}$ and $D_{\Gt_{3}}$ correspond
to the curves with two irreducible components, and $\cmodo{4}$ is
the complement of these three points (see Fig. \ref{moduli}). The
universal family $\umod{4}$ is a del Pezzo surface of degree five
which is obtained by blowing up three points of $\projc \times
\projc$.

(iii) The moduli space $\cmod{5}$ is isomorphic to $\umod{4}$. It
has ten divisors and each of these divisors contains three
codimension two strata. The corresponding $5$-trees are shown in
Fig. \ref{moduli}.

\subsection{Forgetful morphism}
\label{sec_forget}%

We say that $(\GS;p_{1},\cdots,p_{n-1})$ is obtained by forgetting
the labeled point $p_{n}$ of the $n$-pointed stable curve
$(\GS;p_{1},\cdots, p_{n}) $. However, the resulting pointed curve
may well be unstable. This happens when the component $\GS_{v}$ of
$\GS$ supporting $p_{n}$ has only two additional special points. In
this case, we contract this component to its intersection point(s)
with the components adjacent to $\GS_{v}$. With this {\it
stabilization} we extend this map to whole space and obtain
$\pi_{n}: \cmod{n} \to \cmod{n-1}$. There exists a canonical
isomorphism $\cmod{n} \to \umod{n-1}$ commuting with the projections
to $\cmod{n-1}$. In other words, $\pi_{n}: \cmod{n} \to \cmod{n-1}$
can be identified with the universal family of curves. For details,
see \cite{ke,knu}.

\subsection{Automorphisms of $\cmod{n}$}
\label{automorphisms}%

The open stratum $\cmodo{n}$ of the moduli space $\cmod{n}$ can be
identified with the orbit space $((\projc)^n \smin \GD)/PSL_2(\C)$.
The latter orbit space may be viewed as the configuration space of
$(n-3)$ ordered distinct points of  $\projc \smin \{0,1,\infty\}$:
\begin{eqnarray}
\cmodo{n} \cong \{\ve{p} = (z_1,\cdots,z_n) \in \C^{n-3} \mid z_i
\ne z_j \ \forall i \ne j,\ z_{n-2}=0, z_{n-1}= 1, z_n= \infty \},
\nonumber
\end{eqnarray}
where $z_i := [z_i:1]$ are the coordinates of labeled points $p_i$
in an affine chart of $\projc$.

Let $\psi=(\psi_1,\ldots,\psi_{n-3}): \cmodo{n} \to \cmodo{n}$ be a
non-constant holomorphic map. In \cite{kalim}, Kaliman discovered
the following fact:
\begin{thm}[Kaliman, \cite{kalim}]
For $n \geq 4$, every non-constant holomorphic endomorphism
 $\psi=(\psi_1,\cdots,\psi_{n-3})$ of  $M_n$
 is an automorphism and its components $\psi_r$ are of the form
\begin{eqnarray}
\psi_r(\ve{p}) = \frac{z_{\Gr(r)} - z_{\Gr(n-2)}}{z_{\Gr(r)} -
z_{\Gr(n)}} \ \Big/ \ \frac{z_{\Gr(n-2)} - z_{\Gr(n-1)}}{z_{\Gr(n)}
- z_{\Gr(n-1)}}, \ 1 \leq r \leq n-3 \nonumber
\end{eqnarray}
where $\Gr \in S_n$ is a permutation not depending on $r$.
\end{thm}

Kaliman's theorem implies the following corollary.
\begin{thm}[Kaliman \& Lin, \cite{kalim,lin}]
Every holomorphic automorphism of $\cmodo{n}$ is produced by a
certain permutation $\Gr \in S_n$. Hence, $Aut(\cmodo{n}) \cong
S_n$.
\end{thm}

On the other hand,  the permutation group $S_n$ acts on the
compactification $\cmod{n}$ of $\cmodo{n}$ via relabeling: for $\Gr
\in S_n$, there is a holomorphic map $\psi_\Gr$ which is given by
\begin{eqnarray} \label{eqn_auto}
 \psi_\Gr: \curve \mapsto (\GS;\Gr(\ve{p})) := (\GS;p_{\Gr(1)},
 \cdots,p_{\Gr(n)}).
\end{eqnarray}
Therefore, the permutation action given in (\ref{eqn_auto}) forms a
subgroup of holomorphic automorphisms $Aut(\cmod{n})$.

Let $Aut_\sharp(\cmod{n})$ be the group of holomorphic automorphisms
of $\cmod{n}$ that respect the stratification: $\psi \in
Aut_\sharp(\cmod{n})$ maps $D_\Gt$ onto $D_\Gg$ where $\mathrm{dim}
D_{\Gt} = \mathrm{dim} D_{\Gg}$. Kaliman's theorem leads us to the
following immediate corollary.
\begin{thm} \label{cor_auto}
The group $Aut_\sharp(\cmod{n})$ is $S_n$.
\end{thm}

\begin{proof}
Let $\psi \in Aut_\sharp(\cmod{n})$. The restriction of $\psi$  to
the open stratum gives the permutation action on $\cmodo{n}$ since
the automorphism group of the open stratum $\cmodo{n}$ contains only
permutations $\psi_\Gr$. The unicity theorem of holomorphic maps
implies that $\psi = \psi_\Gr$ since they coincide on the open
stratum $\psi|_{\cmodo{n}} = \psi_\Gr |_{\cmodo{n}}$.
\end{proof}

\begin{rem}
Note that the whole group of holomorphic automorphisms
$Aut(\cmod{n})$ is not necessarily isomorphic to $S_n$. For example,
the automorphism group of $\cmod{4}$ is $PSL_2(\C)$. To the best of
our knowledge, there is no systematic exposition of $Aut(\cmod{n})$
for $n >5$.
\end{rem}

\section{Moduli of pointed real curves of genus zero}
\label{ch_real_curves}

A {\it real structure} on complex variety $X$ is an anti-holomorphic
involution $c_X: X\to X$. The fixed point set $\real{X} = \fix(c_X)$
of the involution is called the {\it real part} of the variety (or
of the real structure).

In this section, we introduce the moduli spaces of pointed real
curves of genus zero as the fixed point sets of  real structures on
$\cmod{n}$.

\subsection{Real structures on $\cmod{n}$}
\label{sec_real_curves}%

The moduli space $\cmod{n}$ comes equipped with a natural real
structure. The involution $c: \curve  \mapsto
(\overline{\GS};\mathbf{p})$ acts on $\cmod{n}$. Here a complex
curve $\GS$ is regarded as a pair $\GS=(C,J)$, where $C$ is the
underlying 2-dimensional manifold and $J$ is a complex structure on
it, and $\overline{\GS} = (C,-J)$ is its complex conjugated
pair.\footnote{There is some notational ambiguity: The bar over
$\cmod{n}$ and that over $\overline{\GS}$ refer to different
structures on underlying manifolds: the first one refers to the
compactification of $\cmod{n}$ and second refers to the manifold
with reverse complex structure. Both of these notations are widely
used, we use the bar for both cases. The context should make it
clear which structure is referred to.}

\begin{lem} \label{lem_real_str}
The map $c$ is a real structure on $\cmod{n}$.
\end{lem}

\begin{proof} The differentiability of $c$ follows form the
Kodaira-Spencer construction of infinitesimal deformations. We need
to show that the differential of $c$ is anti-linear at each $\curve
\in \cmod{n}$. It is sufficient to show that it is anti-linear
without taking the quotient with respect to $PSL_2(\C)$.

The infinitesimal  deformations of a nodal curve $\GS$ with divisor
$D_{\mathbf{p}}= p_1+\cdots+p_n$ is canonically identified with the
complex vector space $Ext^{1}_{\cO_{\GS}} (\GO^1_{\GS}
(D_{\mathbf{p}}),\cO_{\GS})$, (see Section \ref{sec_deform}). By
reversing the complex structure on $\GS$, we reverse the complex
structure on the tangent space $Ext^{1}_{\cO_{\GS}} (\GO_{\GS}
(D_{\mathbf{p}}),\cO_{\GS})$ at $\curve$. The differential of the
map $\curve \mapsto (\overline{\GS};\mathbf{p})$
\begin{eqnarray*}
Ext^{1}_{\cO_{\GS}} (\GO_{\GS} (D_{\mathbf{p}}),\cO_{\GS})  \to
Ext^{1}_{\cO_{\overline{\GS}}} (\GO_{\overline{\GS}}
(D_{\mathbf{p}}), \cO_{\overline{\GS}})
\end{eqnarray*}
is clearly anti-linear.
\end{proof}

The subgroup  $Aut_\sharp(\cmod{n})\cong S_n$ of holomorphic
automorphisms acts on $\cmod{n}$ via relabeling as given in
(\ref{eqn_auto}). For each involution $\Gs \in S_{n}$, we have an
additional real structure on $\cmod{n}$:
\begin{equation}
\label{eqn_real_str}%
c_{\Gs}:= c \circ \psi_{\Gs}: \curve \mapsto \curvef.
\end{equation}

\begin{lem}
Every real structure of $\cmod{n}$ preserving the stratification is
of the form  (\ref{eqn_real_str}) where $\Gs \in S_n$ is an
involution.
\end{lem}

\begin{proof}
By their definition, anti-holomorphic automorphisms of $\cmod{n}$
are obtained by composing the principal real structure $c: \curve
\mapsto (\konj{\GS};\ve{p})$ with elements of $Aut(\cmod{n})$. The
real structure $c$ maps each stratum of $\cmod{n}$ onto itself.
Therefore, an anti-holomorphic automorphism $c \circ \psi$ respects
the stratification of $\cmod{n}$ if and only if $\psi \in Aut(\cmod{n})$
respects the stratification of $\cmod{n}$ i.e., each real structure
preserving the stratification of $\cmod{n}$ is given by a certain
involution $\Gs \in S_n$ and is of the form (\ref{eqn_real_str}).
\end{proof}

\subsection{$\Gs$-invariant curves and $\Gs$-equivariant families}

\begin{defn}
\label{def_s_curve}%
An $n$-pointed stable curve $\curve$ is called {\it $\Gs$-invariant}
if it admits a real structure $\conj: \GS \to \GS$ such that
$\conj(p_i) = p_{\Gs(i)}$ for all $i \in \{\Nn\}$.

A family of $n$-pointed stable curves $\pi_S: \cU_S \to S$ is called
{\it $\Gs$-equivariant} if there exist a pair of real structures
\begin{eqnarray*}
\begin{CD}
\cU_S  @> {c_{\cU}} >>   \cU_S \\
@V{\pi_S}VV      @VV{\pi_S}V   \\
  S    @> {c_S}     >>       S.
\end{CD}
\end{eqnarray*}
such that the fibers $\pi^{-1} (s)$ and $\pi^{-1} (c_S(s))$ are
$\GS$ and $\konj{\GS}$ respectively, and $c_\cU$ maps $z \in
\pi^{-1} (s)$ to $z \in \pi^{-1} (c_S(s))$.
\end{defn}

\begin{rem}
If $\curve$ is $\Gs$-invariant, then the real structure $\conj: \GS
\to \GS$ is uniquely determined by the permutation $\Gs$ due to
Lemma \ref{lem_bihol}.
\end{rem}

\begin{lem} \label{lem_real_curve}
If $\pi: \cU_S \to S$ is a $\Gs$-equivariant family, then each
$\curve \in \real{S}$ is $\Gs$-invariant.
\end{lem}

\begin{proof} If $c_{S} (\curve) = (\GS;\ve{p})$ then
there exists a unique bi-holomorphic equivalence $\conj: \GS \to
\overline{\GS}$ (in other words, anti-holomorphic $\conj: \GS \to
\GS$) such that $\conj(p_i) = p_{\Gs(i)}$. The restriction of
$\conj$ on labeled points is an involution. By applying Lemma
\ref{lem_bihol} to $\conj^2$, we determine that $\conj$ is an
involution on $\GS$.
\end{proof}

\subsection{The moduli space of pointed real curves  $\rmod{2k,l}$}
\label{sec_r_mod}

Let $\fix(\Gs)$ be the fixed point set of the action of $\Gs$ on the
labeling set $\{\Nn\}$, and let $\Conj (\Gs)$ be its complement. Let
$|\fix(\Gs)|=l$ and $|\Conj(\Gs)|=2k$.

Let $\Gr \in S_n$, and  $\psi_{\Gr}$ be the corresponding
automorphism of $\cmod{n}$. The conjugation of real structure
$c_\Gs$ with $\psi_{\Gr}$ provides a conjugate real structure
$c_{\Gs'} = \psi_{\Gr} \circ c_{\Gs} \circ \psi_{\Gr^{-1}}$. The
conjugacy  classes of real structures are determined by the
cardinalities $|\fix(\Gs)|$ and $|\Conj(\Gs)|$. Therefore, from now
on, we only  consider $c_{\Gs}$ where
\begin{eqnarray} \label{eqn_sigma}
\Gs = \left(
\begin{array}{ccccccccc}
1   & \cdots & k  & k+1 & \cdots & 2k & 2k+1 & \cdots & 2k+l \\
k+1 & \cdots & 2k &  1  & \cdots & k  & 2k+1 & \cdots & 2k+l \\
\end{array}
\right),
\end{eqnarray}
and $n=2k+l$.

\begin{defn}
For $\Gs$ as above in (\ref{eqn_sigma}), $\Gs$-invariant curves are
called {\it $(2k,l)$-pointed real curves}.

The fixed point set $\fix(c_\Gs)$ is called the {\it moduli space of
$(2k,l)$-pointed real curves} and denoted by  $\rmod{2k,l}$.
\end{defn}

\begin{thm}
\label{thm_r_mod}%
(i) For any $n \geq 3$, $\rmod{2k,l}$ is a smooth real projective
manifold of  dimension $n-3$.

(ii) The universal family of curves $\pi: \umod{n} \to \cmod{n}$ is
a $\Gs$-equivariant family.

(iii) Any  $\Gs$-equivariant family of  $n$-pointed stable curves
over $\pi_S:\cU_S \to S$ is induced by a unique pair of real
morphisms
\begin{eqnarray*}
\begin{CD}
\cU_S  @> \hat{\kappa} >>    \umod{n} \\
@V{\pi_S}VV      @VV{\pi}V   \\
  S    @> \kappa >>          \cmod{n}.
\end{CD}
\end{eqnarray*}

(iv) Let $\fM_\Gs$ be the contravariant functor that sends real
varieties $(B,c_B)$ to the set of $\Gs$-equivariant families of
curves over $B$. The moduli functor $\fM_\Gs$ is represented by the
real variety $(\cmod{n},c_\Gs)$.

(v) Let $\R \fM_\Gs$ be the contravariant functor that sends real
analytic manifolds $R$ to the set of families of $\Gs$-invariant
curves over $R$. The moduli functor $\R\fM_\Gs$ is represented by
the real part $\rmod{2k,l}$ of $(\cmod{n},c_\Gs)$.
\end{thm}

\begin{proof}
(i) The smoothness of the real part of $c_\Gs$ is a consequence of
the implicit function theorem, and $dim_{\R} \ \rmod{2k,l} =
dim_{\C} \ \cmod{n} =n-3$ since the real part $\rmod{2k,l}$ is not
empty.

(ii) The fiber over $\curve \in \cmod{n}$ is $\pi^{-1}(\curve) =
\GS$. We define  real structures on $\cmod{n}$ and $\umod{n}$ as
follows;
\begin{eqnarray*}
\begin{array}{llll}
c_{\Gs}:& \curve &\mapsto& (\konj{\GS};\Gs(\ve{p})),\\
\hat{c}_\Gs :&  z \in \pi^{-1}(\curve)   &\mapsto& z \in
\pi^{-1}((\overline{\GS};\Gs(\ve{p}))),
\end{array}
\end{eqnarray*}
The real structures $c_\Gs, \hat{c}_\Gs$ satisfy the conditions of
$\Gs$-equivariant families in Definition \ref{def_s_curve}.

(iii) Due to Knudsen's theorem (see Section \ref{stratacomplex}),
each of the morphisms $\kappa: S \to \cmod{n}$ and $\hat{\kappa}:
\cU_S \to \umod{n}$ are unique. Therefore, they are the same as
$c_{\Gs} \circ \kappa \circ c_{S}: S \to \cmod{n}$ and $\hat{c}_\Gs
\circ \hat{\kappa} \circ c_\cU: \cU_S \to \umod{n}$. Hence, the
morphisms $\kappa, \hat{\kappa}$ are real.

(iv)-(v) The statements follows from (iii) and the definition of
moduli functors.
\end{proof}

It was believed that the real locus $\rmod{2k,l}$ does not represent 
any moduli functor for $k \ne 0$, but it has only a meaning in
operadic setting. Theorem \ref{thm_r_mod} shows the contrary.

\section{Stratification of $\rmod{2k,l}$}
\label{sec_stratification_sp}%

A stratification for $\rmod{2k,l}$ can be obtained by the
stratification of $\cmod{n}$ given in Section \ref{stratacomplex}.

\begin{lem} \label{lem_inv_tree}
Let $\Gg$ and $\overline{\Gg}$ be the dual trees of $\curve$ and
$\curvef$ respectively.

(i) If $\Gg$ and $\overline{\Gg}$ are not isomorphic, then the
restriction of $c_{\Gs}$ on the union of complex strata $D_{\Gg}
\bigcup D_{\overline{\Gg}}$ gives a real structure with empty real
part.

(ii) If  $\Gg$ and $\overline{\Gg}$ are isomorphic, then the
restriction of  $c_{\Gs}$ on $D_{\Gg}$ gives a real structure whose
corresponding real part $\real{D_{\Gg}}$ is the intersection of
$\rmod{2k,l}$ with $D_{\Gg}$.
\end{lem}

\begin{proof}
(i) Since $\Gg$ and $\overline{\Gg}$ are not isomorphic, $D_{\Gg}$
and $D_{\overline{\Gg}}$ are disjoint complex strata. The
restriction of $c_{\Gs}$ on $D_{\Gg} \bigcup D_{\overline{\Gg}}$
swaps the strata. Therefore, the real part of this real structure is
empty.

(ii) Since $\Gg$ and $\overline{\Gg}$ are isomorphic, the
$n$-pointed curves $\curve$ and $\curvef$ are in the same stratum
$D_{\Gg}$. Therefore, the restriction of $c_{\Gs}$ on $D_{\Gg}$ is a
real structure. The real part $\real{D_\Gg}$ of the
$\Gs$-equivariant family $D_\Gg$ is $\rmod{2k,l} \bigcap D_{\Gg}$
since $\rmod{2k,l}= \fix(c_\sigma)$.
\end{proof}

\begin{defn}
A tree $\Gg$ is called {\it $\Gs$-invariant} if it is isomorphic to
$\bar{\Gg}$.  We denote  the set of $\Gs$-invariant $n$-trees  by
$\tree(\Gs)$.
\end{defn}

\begin{thm} \label{cor_strata}
$\rmod{2k,l}$ is stratified by real analytic subsets
$\real{D_{\Gg}}$ where $\Gg \in \tree(\Gs)$.
\end{thm}

Although the notion of $\Gs$-invariant trees leads us to a
combinatorial stratification of $\rmod{2k,l}$ as given in Theorem
\ref{cor_strata}, it does not give a stratification in terms of
connected strata. For a $\Gs$-invariant $\Gg$, the real part of the
stratum $\R D_{\Gg}$ has many connected components. In the next
subsection, we refine this stratification by using the spaces of
$\Z_2$-equivariant point configurations in the projective line
$\projc$.

\subsection{Spaces of $\Z_2$-equivariant point configurations
in $\projc$} \label{sec_conf_space}

Let $z := [z:1]$ be the affine coordinate on $\projc$. Consider the
upper half-plane $\H^+ = \{z \in \projc \mid \Im{(z)} >0 \}$ (resp.
lower half plane  $\H^- = \{z \in \projc \mid \Im{(z)} <0 \}$) as a
half of the $\projc$ with respect to $z \mapsto \bar{z}$, and the
real part $\projr$ as its boundary. Denote by $\H$ the compact disc
$\H^{+} \cup \projr$.

\subsubsection{Irreducible $(2k,l)$-pointed real curves and
configuration spaces} \label{sec_conf_ir}%
Let $\curve$ be an irreducible $(2k,l)$-pointed real curve. As a
real curve, $\GS$ is isomorphic to $\projc$ with real structure
which is either $z \mapsto \bar{z}$ or $z \mapsto -1/\bar{z}$ (see,
for example \cite{m2}). The real structure $z \mapsto -1/\bar{z}$
has empty real part. Notice that a $(2k,l)$-pointed curve with empty
real part is possible only when $\fix(\Gs) = \emp$ i.e., $l=0$.

We will consider the  spaces of real curves with non-empty and empty
real parts as separate cases.

\paragraph*{\bf Case I. {\em Configurations in $\projc$ with non-empty
real part.}} %
Each finite subset $\ve{p}$ of $\projc$ which is
invariant under the real structure $z \mapsto \bar{z}$ inherits
additional structures:

\begin{itemize}
\item[{\bf I.}] {\it An oriented cyclic ordering on $\fix(\Gs)$}:
For any point $p \in (\ve{p}\bigcap \projr)$ there is unique $q  \in
(\ve{p} \bigcap \projr)$ which follows the point $p$ in the positive
direction of $\projr$ (the direction in which the coordinate
$x:=[x:1]$ on $\projr$ increases).

\noindent The elements of $\ve{p}$ are labeled, therefore the cyclic
ordering can be seen as a linear ordering on $(\ve{p} \bigcap
\projr) \smin \{p_n\}$. This linear ordering gives an {\it oriented
cyclic ordering} on $\fix(\Gs) =\{2k+1,\cdots,n\}$ which we denote
by $\{i_1\} < \cdots < \{i_{l-1}\} <\{i_l\}$ where $i_l=n$.

\item[{\bf II.}] {\it A  2-partition on $\Conj(\Gs)$}: The subset
$\ve{p} \bigcap (\projc \smin \projr)$ of $\ve{p}$ admits a
partition into two disjoint subsets $\{p_i \in \H^\pm\}$.  This
partition gives an {\it ordered 2-partition} $\Conj^\pm := \{ i \mid
p_i \in \H^\pm\}$ of $\Conj(\Gs)$. The subsets $\Conj^\pm$ are
swapped by the permutation $\Gs$.
\end{itemize}

The set of data 
\begin{equation*}
o:=\{(\projc,z \mapsto \bar{z});\Conj^\pm;
\fix(\Gs) = \{\{i_1\} < \cdots < \{i_{l}\}\}\} 
\end{equation*}
is called the {\it oriented combinatorial type} of the $\Z_2$-equivariant 
point configuration $\ve{p}$ on $(\projc,z \mapsto \bar{z})$.

The oriented combinatorial types of equivariant point configurations
on $(\projc,z \mapsto \bar{z})$ enumerate the connected components
of the space $\cspp{2k,l}$ of $k$ distinct pairs of conjugate points
on $\H^+ \bigcup \H^-$ and $l$ distinct points on $\projr$:
\begin{eqnarray*}
\cspp{2k,l} &:=& \{ (p_1,\cdots,p_{2k};q_{2k+1},\cdots,q_{2k+l})
\mid
p_i \in \projc \smin \projr, \ p_i = p_j \Leftrightarrow i=j, \\
&&  p_i = \bar{p}_{j} \Leftrightarrow i= \Gs(j) \  \& \ q_i \in
\projr, \ q_i =  q_j  \Leftrightarrow i=j \}.
\end{eqnarray*}

The number of connected components of $\cspp{2k,l}$ is $2^{k}
(l-1)!$.\footnote{ Here we use the convention $n!=1$ whenever $n
\leq 0$.} They are all pairwise diffeomorphic; natural
diffeomorphisms are given by $\Gs$-invariant relabeling.

Let $z := [z:1]$ be the affine coordinate on $\projc$ and $x :=
[x:1]$ be affine coordinate on $\projr$. The action of  $SL_2(\R)$
on $\H$ is given by
\begin{equation*}
SL_2(\R) \times \H \to \H,\ \ (\Lambda,z) \mapsto \Lambda(z) =
\frac{az+b}{cz+d},\ \ \Lambda =  \left(
\begin{array}{cc}
a   &   b \\
c   &   d \\
\end{array}
\right) \in SL_2(\R)
\end{equation*}
in affine coordinates. It induces an isomorphism $SL_2(\R)/\pm I \to
Aut(\H)$. The  automorphism group $Aut(\H)$ acts on $\cspp{2k,l}$
\begin{equation*}
\Lambda: (z_1,\cdots,z_{2k};x_{2k+1},\cdots,x_{2k+l}) \mapsto
(\Lambda(z_1),\cdots,\Lambda(z_{2k});\Lambda(x_{2k+1}),\cdots,
\Lambda(x_{2k+l})).
\end{equation*}
It preserves each of the connected components of $\cspp{2k,l}$. This
action is free when $2k+l \geq 3$, and it commutes with
diffeomorphisms given by $\Gs$-invariant relabellings. Therefore, the
quotient space $\cspq{2k,l}:=\cspp{2k,l}  / Aut(\H)$ is a manifold
of dimension $2k+l-3$ whose connected components are pairwise
diffeomorphic.

In addition to the automorphisms considered above, there is a
diffeomorphism $-\I$ of $\cspp{2k,l}$ which is given in affine
coordinates as follows.
\begin{equation} \label{eqn_-1}
-\I: (z_1,\cdots,z_{2k};x_{2k+1},\cdots,x_{2k+l}) \mapsto
(-z_1,\cdots,-z_{2k};-x_{2k+1},\cdots,-x_{2k+l}).
\end{equation}
Consider the quotient space $\csp{2k,l} = \cspp{2k,l}/(-\I)$. Note
that, $-\I$ interchanges components with {\it reverse} combinatorial
types. Namely, the combinatorial type $\konj{o}$ of $-\I(\ve{p})$ is
obtained from the combinatorial type $o$ of $\ve{p}$ by reversing
the cyclic ordering on $\fix(\Gs)$ and swapping $\Conj^+$ and
$\Conj^-$. The equivalence classes of oriented combinatorial types
with respect to the action of $-\I$ are called {\it un-oriented
combinatorial types} of $\Z_2$-equivariant point configurations on
$(\projc,z \mapsto \konj{z})$. The un-oriented combinatorial types
enumerate the connected components of $\csp{2k,l}$.

The diffeomorphism $-\I$ commutes  with each $\Gs$-invariant
relabeling and normalizing action of $Aut(\H)$. Therefore, the
quotient space $\csq{2k,l} := \csp{2k,l} /Aut(\H)$ is a manifold of
dimension $2k+l-3$, its connected components are diffeomorphic to
the components of $\cspq{2k,l}$, and, moreover, the quotient map
$\cspq{2k,l} \to \csq{2k,l}$ is a trivial double covering.

\paragraph*{\bf Case II. {\em Configurations in $\projc$ with empty
real part.}}%
Let $\curve$ be an irreducible $(2k,0)$-pointed real
curve and let $\real{\GS} = \emp$. Such a pointed real curve is
isomorphic to $(\projc,\ve{p})$ with real structure $\conj: z
\mapsto -1/\bar{z}$.

The group of automorphisms of $\projc$ which commutes with $\conj$
is
\begin{eqnarray}
Aut(\projc,\conj)\cong SU(2) := \left\lbrace  \left(
\begin{array}{cc}
a   &   b \\
-\bar{b}   &   \bar{a} \\
\end{array}
\right) \in SL_2(\C) \right\rbrace.  \nonumber
\end{eqnarray}
Thus, the group $Aut(\projc,\conj)$ acts naturally on the space
\begin{equation*}
Conf_{(2k,0)}^\emp :=\{(z_1,\cdots,z_{2k}) \mid z_i =
-1/\bar{z}_{i+k}\}
\end{equation*}
of $\Z_2$-equivariant point configurations on $(\projc,z \mapsto
-1/\bar{z})$. For $k \geq 2$, the action is free and the quotient
$\dsq{2k,0} := Conf_{(2k,0)}^\emp / Aut(\projc,\conj)$ is a $2k-3$
dimensional connected manifold.

The {\it combinatorial type}  of $\Z_2$-equivariant point
configurations on $(\projc,z \mapsto -1/\bar{z})$ is unique and
given by the topological type of the real structure $z \mapsto
-1/\bar{z}$.

\subsubsection{A normal position of $\Z_2$-equivariant point
configurations on $\projc$} \label{sec_normal_pos} By using the
automorphisms we can make the following choices for the
representatives of the points in $\cspq{2k,l}$ and $\dsq{2k,0}$.

\paragraph*{\bf Case I. {\em Configurations in $\projc$ with non-empty
real part.}}%
Every element in $\cspq{2k,l}$ is represented by $(\projc,\ve{p})$
with $\ve{p} \in \cspp{2k,l}$. In order to calibrate the choice by
$Aut(\H)$, consider an isomorphism $(\projc,\ve{p}) \mapsto
(\projc,\ve{p}')$ which puts the labeled points in the following
normal position $\ve{p}' \in \projc$.

\begin{itemize}
\item[\bf{(A)}]
In the case $l \geq 3$, the three consecutive labeled points
$(p'_{i_{l-1}},p'_{n},p'_{i_{1}})$ in $\projr$ are put in the
position $x'_{i_{l-1}}= 1,x'_{n}=\infty,x'_{i_{1}}=0$. We then
obtain
\begin{eqnarray*}
\ve{p}' = (z_1,\cdots,z_k,\overline{z}_{1},\cdots,\overline{z}_{k},
x_{2k+1},\cdots,x_{2k+l-1},\infty).
\end{eqnarray*}

\item[\bf{(B)}] In the case $l =1,2$, the three labeled points
$\{p_k,p_{2k},p_n\}$ are put in the position  $\{\pm\im,\ \infty\}$.
Then,
\begin{eqnarray*}
\ve{p}' = \left\{
\begin{array}{ll}
(z_1,\cdots,z_{k-1},\epsilon
\im,\overline{z}_{1},\cdots,\overline{z}_{k-1},
-\epsilon \im, x_{2k+1},\infty) & \mathrm{if} \ l=2, \\
(z_1,\cdots,z_{k-1}, \epsilon
\im,\overline{z}_{1},\cdots,\overline{z}_{k-1}, -\epsilon
\im,\infty) & \mathrm{if} \ l=1
\end{array} \right.
\end{eqnarray*}
where $\epsilon =\pm$.

\item[\bf{(C)}] In the case $ l =0$, the labeled points
$\{p_k,p_{2k}\}$  are fixed at  $\{ \pm \im \}$ and  $p_{i}$ where
$\{i\} = \{k-1,2k-1\} \bigcap \Conj^+$ is placed on the interval
$]0,\im[ \subset \H^+$. Then,
\begin{eqnarray*}
\mathbf{p}' = (z_1,\cdots,z_{k-2}, \epsilon_1 \Gl \im, \epsilon_2
\im, \overline{z}_{1},\cdots,\overline{z}_{k-2}, -\epsilon_1 \Gl
\im, -\epsilon_2 \im).
\end{eqnarray*}
where  $\Gl \in ]0,1[$ and $\epsilon_i=\pm, i=1,2$.
\end{itemize}

\begin{rem}[{\bf A}$'$]
In case of $k>1$ and $l>3$, we can consider the alternative map
which puts the labeled points $\ve{p}'$ into the following normal
positions. In this case, three labeled points $\{p_k,p_{2k},p_n\}$
can be put at  $\{z_k,z_{2k},x_n\}=\{\pm\im,\infty\}$ by the action
of $Aut(\H)$. We then obtain
\begin{eqnarray*}
\mathbf{p}' = (z_1,\cdots,z_{k-1},\epsilon
\im,\overline{z}_{1},\cdots, \overline{z}_{k-1}, -\epsilon \im,
x_{2k+1},\cdots,x_{2k+l-1}, \infty)
\end{eqnarray*}
where $\epsilon=\pm$.
\end{rem}

\paragraph*{\bf Case II. {\em Configurations in $\projc$ with empty real
part.}} %
Let $k \geq 2$.
Every element of $\dsq{2k,0}$  is represented by
$(\projc,\ve{p})$ with $\ve{p} \in Conf^\emp_{(2k,0)}$. Calibrating
the choice by $Aut(\projc,\conj)$, consider an isomorphism
$(\projc,\ve{p}) \mapsto (\projc,\ve{p}')$ which puts the labeled
points of $\curve$ in the following normal position $\ve{p}' \in
\projc$.
\begin{itemize}
\item[{\bf (D)}]
\begin{eqnarray*}
\ve{p} = (z_1,\cdots,z_{k-2}, \Gl \im ,  \im,
\frac{-1}{\overline{z}_{1}},\cdots,\frac{-1}{\overline{z}_{k-2}}, -
\frac{\im}{\Gl},- \im)
\end{eqnarray*}
where $\Gl \in ]-1,1[$.
\end{itemize}

\subsection{O/U-planar trees:  one-vertex case.}
\label{def_o-pl} An {\it oriented planar (o-planar) structure } on
the one-vertex $n$-tree $\Gt$ is one of the two possible sets of
data
\begin{eqnarray*}
o&:=&  \left\{
    \begin{array}{ll}
     \{ \real{\GS} \ne \emp; \mathrm{a \ \Gs-equivariant\ two-partition\
        \Conj^\pm \ of }\ \Conj(\Gs); \\
      \mathrm{an\ oriented  \ cyclic \ ordering \  on}\  \fix(\Gs)\}, \\
      \{ \real{\GS}  = \emp\}.
    \end{array}
        \right. \nonumber
\end{eqnarray*}
We denote the o-planar trees by $(\Gt,o)$.

An {\it un-oriented planar (u-planar) structure} $u$ on the
one-vertex $n$-tree $\Gt$ is a pair of reverse o-planar structures
$\{o,\konj{o}\}$ when $\real{\GS} \ne \emp$, and equal to the
o-planar structure when $o=\{ \real{\GS} = \emp\}$. We denote the
u-planar trees by $(\Gt,u)$.

\subsubsection{O/U-planar trees and connected components of configuration
spaces}%
As shown in Section \ref{sec_conf_ir}, each connected component of
$\csq{2k,l}$ for $l>0$ (resp. $\csq{2k,0} \bigcup \dsq{2k,0}$ for
$l=0$) is associated to a unique u-planar tree since  the un-oriented 
combinatorial types of $\Z_2$-equivariant point configurations are 
encoded by the same set of data. We denote the connected components 
of $\csq{2k,l}$ (and $\csq{2k,0} \bigcup \dsq{2k,0}$) by $\csq{\Gt,u}$. 
Similarly, each connected component of $\cspq{2k,l}$ is associated to 
a unique o-planar tree. We denote the connected components of 
$\cspq{2k,l}$ by $\csq{\Gt,o}$.

\subsubsection{Connected components of $\rmodo{2k,l}$}
Every $\Z_2$-equivariant point configuration defines a
$(2k,l)$-pointed real curve. Hence, we define
\begin{equation} \label{eqn_dif}
\Xi: \bigsqcup_{(\Gt,u)} \csq{\Gt,u} \to \rmodo{2k,l}
\end{equation}
which maps  $\Z_2$-equivariant point configurations to the
corresponding isomorphism classes of irreducible $(2k,l)$-pointed
curves.

\begin{lem}
\label{lem_conf_sp} %
(i) The map $\Xi$ is a diffeomorphism.

(ii) The configuration space $\csq{\Gt,u}$ is diffeomorphic to
\begin{itemize}
\item $((\H^+)^k \smin \GD) \times \square^{l-3}$ when $l >2$,
\item $((\H^+ \smin \{\im\})^{k-1}  \smin \GD )\times \square^{l-1}$ when
$l=1,2$,
\item $((\H^+ \smin \{\im, \im/2 \})^{k-2}  \smin \GD )\times
\square^{1}$ when $l=0$ and the type of real structure is $(\projc,
z \mapsto \bar{z})$,
\item $((\projc \smin \{\im,\im/2,-\im/2,-\im\})^{k-2}  \smin (\GD\cup\GD^c))
\times \square^1$ when $l=0$ and the type of the real structure is
$(\projc, z \mapsto -1/\bar{z})$.
\end{itemize}
Here, $\GD$ is the  union of all diagonals $z_i\ne z_j (i\ne j)$,
$\GD^c$ is the union of all cross-diagonals $z_i\ne -\frac{1}{\bar
z_j} (i\ne j)$, and $\square^l$ is the $l$-dimensional open simplex.
\end{lem}

\begin{proof}
(i) The map $\Xi$ is clearly smooth. It is surjective since any
$(2k,l)$-pointed irreducible curve is isomorphic either to
$(\projc,z\mapsto \bar{z})$ or $(\projc,z \mapsto -1/\bar{z})$ with
a $\Z_2$-equivariant point configuration $\ve{p}$ on it. It is
injective since the group of holomorphic automorphisms commuting
with the real structure $z \mapsto \bar{z}$ is generated by
$Aut(\H)$ and $-\I$, and the group of holomorphic automorphisms
commuting with the real structure $z \mapsto -1/\bar{z}$ is
$Aut(\projc,\conj)$. These automorphisms are taken into account
during construction of the configuration spaces.

(ii) As it is shown in Section \ref{sec_conf_ir}, $\csq{\Gt,u}$ is
the quotient $\csq{\Gt,o} \bigsqcup \csq{\Gt,\bar{o}}/(-\I$). The
spaces $\csq{\Gt,u}$ and $\csq{\Gt,o}$ are clearly diffeomorphic. To
replace $\csq{\Gt,u}$ by $\csq{\Gt,o}$, we  choose an o-planar
representative for each one-vertex u-planar tree $(\Gt,u)$ among
$(\Gt,o), (\Gt,\bar{o})$ as follows.
\begin{itemize}
\item $l \geq 3$ case: Let $(\Gt,o)$ be the representative of
$(\Gt,u)$ for which $p_{2k+1} < p_{n-1} < p_n$ with respect  to the
cyclic ordering on $\fix(\Gs)$.
\item $l = 0,1,2$ case: Let $(\Gt,o)$ be the representative of
$(\Gt,u)$ for which $k \in \Conj^+$.
\end{itemize}

We put the $\Z_2$-equivariant point configurations into a normal
position as in \ref{sec_normal_pos}. The parameterizations stated in
(ii) for $l>0$ cases follow from ({\bf A}) and ({\bf B}) of
\ref{sec_normal_pos}.  In the case of $l=0$ and $(\projc,z \mapsto
\bar{z})$, according to ({\bf C}) the configuration space
$\csq{\Gt,o}$ is a locally trivial fibration  over $\square^{1} =
]0,1[$ whose fibers over $\Gl \in \square^{1}$ are $(\H^+ \smin
\{\im,\Gl \im\})^{k-2} \smin \GD$. Similarly, in the case of $l=0$
and $(\projc,z \mapsto -1/\bar{z})$, according to ({\bf D}) the
configuration space $\csq{\Gt,u}$ is a locally trivial fibration
over $\square^{1} = ]-1,1[$ whose fibers over $\Gl \in \square^{1}$
are $(\projc \smin \{\im,\Gl \im, -\im/\Gl,-\im\})^{k-2} \smin \GD$.
Since the bases of these locally trivial fibrations are
contractible, they are trivial fibrations, and the result follows.
\end{proof}

\subsection{Pointed real curves and o-planar trees}
This section extends the notions of o/u-planar structures to the all
$\Gs$-invariant trees.

\subsubsection{Notations}
\label{sec_not1}%
Let $\curve \in \real{D_\Gg}$ for some  $\Gg \in \tree(\Gs)$ and
$\conj: \GS \to \GS$ be the real structure on $\GS$. We denote the
set of real components $\{v \in V_\Gg \mid \conj(\GS_v) = \GS_v \}$
of $\GS$ by $V_{\Gg}^\R$. If this set is empty, then $\R \GS$ is an
isolated real node; we call the edge of $\Gg$ representing the
isolated real node the {\it special invariant edge}.

Two vertices $v,\bar{v} \in V_{\Gg} \smin V_{\Gg}^\R$ are said to be
{\it conjugate} if the real structure $\conj$ maps the components
$\GS_v$ and $\GS_{\bar{v}}$ onto each other. Similarly, we call the
flags $f,\bar{f} \in F_{\Gg} \smin \{f_* \in F_{\Gg} \mid p_{f_*}
\in \real{\GS}\}$ conjugate if $\conj$ swaps the corresponding
special points $p_f,p_{\bar{f}}$.

\subsubsection{ O/U-planar trees: general case.}
\label{sec_o-planar}%
Let $\Gg \in \tree(\Gs)$. Each $(2k,l)$-pointed real curve $\curve
\in \real{D_\Gg}$ with $\real{\GS} \ne \emp$ admits additional
structures:
\begin{itemize}
\item[{\bf I.}] If $V_\Gg^\R \ne \emp$, one can fix an
oriented combinatorial type of point configurations on each real
component $\GS_v$. Namely, for $v \in V_\Gg^\R$ with $\real{\GS} \ne
\emp$, an oriented combinatorial type $o_v$ is given by an oriented
cyclic ordering on the set of {\it invariant flags} $F_\Gg^\R(v)=\{f
\mid p_f \in \real{\GS_v}\}$ and an  ordered $2$-partition
$F_{\Gg}^\pm(v) = \{i \mid p_i \in \H^\pm\}$ of $F_\Gg(v) \smin
F_\Gg^\R(v)$. In the case of $V_\Gg^\R = \{v\}$ and $\real{\GS} =
\emp$, the o-planar structure is simply given by the type of the
real structure, i.e., $\real{\GS_{v}}=\emp$ (as in \ref{sec_conf_ir}).

\item[{\bf II.}] If $V_\Gg^\R = \emp$, then one can fix an
ordering of the flags of the special invariant edge: $f_e \mapsto
\pm$, $f^e \mapsto \mp$.
\end{itemize}
This additional structures motivate the following definition.

\begin{defn}
An {\it o-planar} structure on $\Gg \in \tree(\Gs)$ is the following
set of data,
\begin{eqnarray}
o:= \left\{
    \begin{array}{ll}
      \{ (\Gg_v,o_v) \mid v \in V_{\Gg}^\R \} & \mathrm{when}\ V_\Gg^\R \ne \emp , \\
      \{ \mathrm{ordering}\ \{f_e,f^e\} \to \{\pm\} \}, & \mathrm{when\ }
      V_{\Gg}^\R = \emp,\\
   \end{array}
        \right.
    \nonumber
\end{eqnarray}
where $(\Gg_v,o_v)$ is an o-planar structure on the one-vertex tree
$\Gg_v$ for  $v \in V_{\Gg}^\R$, and $e=(f_e,f^e)$ is the special
invariant edge of $\Gg$ when $V_\Gg^\R =\emp$.

Similarly, an {\it u-planar} structure on $\Gg$ is the following set
of data,
\begin{eqnarray}
u:= \left\{
    \begin{array}{ll}
          \{ (\Gg_v,u_v) \mid v \in V_{\Gg}^\R \} &
                \mathrm{when}\ V_\Gg^\R \ne \emp, \\
          \{\emp\} & \mathrm{when}\ V_\Gt^\R = \emp,
    \end{array}
    \right.
    \nonumber
\end{eqnarray}
where $(\Gg_v,u_v)$ is an u-planar structure on the one-vertex tree
$\Gg_v$ for  $v \in V_{\Gg}^\R$. We denote u-planar trees by
$(\Gg,u)$.
\end{defn}

\subsubsection{Notations.}
We associate the subsets of vertices $V_{\Gg}^\pm$ and flags
$F_\Gg^\pm$ to every o-planar tree $(\Gg,o)$ when $V_\Gg^\R \ne
\emp$ as follows. Let $v_1 \in V_\Gg \smin V_\Gg^\R$ and let $v_2
\in V_\Gg^\R$ be the closest invariant vertex to $v_1$ in $||\Gg||$.
Let $f \in F_{\Gg}(v_2)$ be in the shortest path connecting the
vertices $v_1$ and $v_2$. The set $V_\Gg^\pm$ is the subset of
vertices $v_1 \in V_\Gg \smin V_\Gg^\R$ such that the flag $f$
(defined as above) is in $F_{\Gg}^\pm(v_2)$. The subset of flags
$F_\Gg^\pm$ is defined as $\dd_\Gg^{-1}(V_{\Gg}^\pm)$.

Similarly, we associate the subsets of vertices $V_{\Gg}^\pm$ and
flags $F_\Gg^\pm$ to every o-planar tree $(\Gg,o)$ when $V_\Gg^\R =
\emp$ : Let $f_e \mapsto +$, $\dd_\Gg(f_e)= v_e$ and $\dd_\Gg(f^e)=
v^e$. The set $V_\Gg^+$ is the subset of vertices $v$ that are
closer to $v_e$ than $v^e$ in $||\Gt||$. The complement $V_\Gg \smin
V_\Gg^+$ is denoted by $V_\Gg^-$. The subset of flags $F_\Gg^\pm$ is
defined as $\dd_\Gg^{-1}(V_{\Gg}^\pm)$.

\subsection{O-planar trees and their configuration spaces}
We associate a product of configuration spaces of $\Z_2$-equivariant
point configurations $\csq{\Gt_v,o_v}$ and moduli space of pointed
complex curves $\cmod{|v|}$ to each o-planar tree $(\Gt,o)$:
\begin{eqnarray*}
\csq{\Gt,o} &:=&
     \left\lbrace
      \begin{array}{ll}
\prod_{v \in V_{\Gt}^\R} \csq{\Gt_v,o_v} \times \prod_{v \in
V_\Gt^+} \cmodo{|v|} &  \mathrm{when}\
V_\Gt^\R \ne \emp \ \mathrm{and}\  \real{\GS} \ne \emp, \\
 \csq{\Gt_{v_r},o_{v_r}} \times
\prod_{\{v,\bar{v}\} \subset V_\Gt \smin V_{\Gt}^\R} \cmodo{|v|} &
\mathrm{when}\
V_\Gt^\R \ne \emp \ \mathrm{and}\  \real{\GS}  =  \emp, \\
\prod_{v \in  V_\Gt^+} \cmodo{|v|} & \mathrm{when }\ V_\Gt^\R =\emp. \\
      \end{array}
     \right.  \\
\end{eqnarray*}
For the case $\real{\GS} = \emp$ , $v_r$ is the vertex corresponding
to the unique real component of point curves, and the product runs
over the un-ordered pairs of conjugate vertices belonging to $V_\Gt
\smin V_\Gt^\R$ i.e., $\{v,\bar{v}\} = \{\bar{v},v\}$.

For each u-planar tree $(\Gt,u)$, we first choose an o-planar
representative and then put $\csq{\Gt,u}=\csq{\Gt,o}$. Note that the
so defined space $\csq{\Gt,u}$ does not depend on the o-planar
representatives up to an isomorphism.

\begin{lem} \label{lem_conf_sp3}
Let $\Gg \in \tree(\Gs)$. The real part $\real{D_\Gg}$ is
diffeomorphic to $\bigsqcup_{(\Gg,u)} \csq{\Gg,u}$ where the
disjoint union is taken over all possible u-planar structures of
$\Gg$.
\end{lem}

\begin{proof}
The complex stratum $D_\Gg$ is diffeomorphic to the product $\prod_{v
\in V_{\Gg}} \cmod{|v|}$. The real structure $c_\Gs: D_\Gg \to
D_\Gg$ maps the factor $\cmod{|v|}$ onto $\cmod{|\bar{v}|}$ when $v$
and $\bar{v}$ are conjugate vertices, and maps the factor
$\cmod{|v|}$ onto itself when $v \in V_\Gg^\R$. Therefore, the real
part $\real{D_\Gg}$ of $c_\Gs$ is given by
\begin{eqnarray*}
\prod_{v \in V_{\Gg}^\R} \csq{2k_v,l_v} \times \prod_{\{v,\bar{v}\}
\subset V_\Gg \smin V_{\Gg}^\R} \cmod{|v|}&
\mathrm{when} \ |V_{\Gg}^\R| >1,  \\
(\csq{|v_r|,0} \bigsqcup B_{(|v_r|,0)}) \times \prod_{\{v,\bar{v}\}
\subset V_\Gg \smin V_{\Gg}^\R} \cmod{|v|}&
\mathrm{when} |V_\Gg^\R|=1,  \\
 \prod_{\{v,\bar{v}\} \subset V_\Gg \smin V_{\Gg}^\R} \cmod{|v|}&
 \mathrm{when} |V_\Gg^\R|=0,
\end{eqnarray*}
where $k_v = |F^+_{\Gg}(v)|$ and $l_v = |F_{\Gg}^\R(v)|$. The
decompositions of the spaces $\csq{2k,l}$ and $\csq{2k,0} \bigsqcup
B_{(2k,0)}$ into their connected components are given in Lemma
\ref{lem_conf_sp}.
\end{proof}

\begin{thm} \label{thm_real_strata}
$\rmod{2k,l}$ is stratified by $\csq{\Gg,u}$.
\end{thm}

\begin{proof}
The moduli space $\rmod{2k,l}$ can be stratified by $\real{D_\Gg}$
due to Theorem \ref{cor_strata}. The claim directly follows from the
decompositions of open strata of $\rmod{2k,l}$ into their connected
components given in Lemma \ref{lem_conf_sp3}.
\end{proof}

\subsection{Boundaries of the strata}
\label{sec_stratification}%
In this section we investigate the adjacency of the strata
$\csq{\Gg,u}$ in $\rmod{2k,l}$. We start with considering the
complex situation in order to introduce natural coordinates near
each codimension one stratum $\csq{\Gg,u}$.

\subsubsection{Intermezzo: Coordinates  around the codimension
one strata} \label{intermezzo} Let $\Gg$ be a 2-vertex $n$-tree
given by  $V_\Gg = \{v_e,v^e\}$, $F_\Gg(v^e) =
\{i_1,\cdots,i_s,f^e\}$ and $F_\Gg(v_e) = \{f_e,i_{s+1},
\cdots,i_{n-1}, n\}$. Let $(z,w):=[z:1] \times [w:1]$ be affine
coordinates on $\projc \times \projc$. Here we introduce coordinates
around $D_\Gg$.

Consider a neighborhood $V \subset D_\Gg$ of a nodal $n$-pointed
curve $(\GS^o,\ve{p}^o) \in D_\Gg$. Any $\curve \in V$ can be
identified with a nodal curve $\{(z -z_{f_e})\cdot(w-w_{f^e})=0\}$
in $\projc \times \projc$ with special points $\ve{p}_{v_e} =
(a_{f_e},a_{i_{s+1}},\cdots,a_{i_{n-1}},a_n) \subset
\{w-w_{f^e}=0\}$ and $\ve{p}_{v^e} =  (b_{f^e},
b_{i_1},\cdots,b_{i_s}) \subset \{z-z_{f_e} =0 \}$. In order to
determine a nodal curve in $\projc \times \projc$  and the position
of its special points uniquely defined by $\curve$, we make the
following choice. Firstly, we fix three labeled points $a_{i_{s+1}},
a_{i_{n-1}}, a_n$ on the line $\{w-w_{f^e}=0\}$ whenever $|v_e|>3$,
and three special points $a_{f_e}, a_{i_{n-1}}, a_n$ whenever
$|v_e|=3$. Secondly, we fix three special points $b_{f^e},
b_{i_{1}}, b_{i_{s}}$ on $\{z-z_{f_e}=0\}$. Finally, we choose
$a_{i_{s+1}}=(0,0)$, $a_{i_{n-1}}=(1,0)$, $a_n=(\infty,0)$ for
$|v_e|>3$; $a_{f_e}=(0,0)$, $a_{i_{n-1}}=(1,0)$, $a_n=(\infty,0)$
for $|v_e|=3$; and $b_{i_1}=(z_{f_e},1)$,
$b_{i_s}=(z_{f_e},\infty)$, $b_{f^e} = (z_{f_e},0)$. Then the
components $z$ and $w$ of the special points provide a coordinate
system in $V$; in particular, for $|v_e|>3$ such a coordinate system
is formed by $z_{f^e}$, $z_{i_*}$ with $i_*=i_{s+2},\cdots,i_{n-2}$,
and  $w_{j_*}$ with $j_*=i_2,\cdots,i_{s-1}$.

We now consider a family of $n$-pointed curves over $V$ times the
$\epsilon$-ball $B_\epsilon = \{|t| < \epsilon\}$. It is given by a
family curves $\{ (z-z_{f_e}) \cdot w  +t=0 \mid t \in B_\epsilon
\}$ in $\projc \times \projc$. The labeled points $(z_i,w_i),
i=1,\cdots n$ on these curves are chosen in the following way. If
$|v_e|>3$, we put $(z_{i_{1}},w_{i_{1}})=(z_{f_e}-t,1)$,
$(z_{i_s},w_{i_s}) = (z_{f_e},\infty)$, $(z_{i_{s+1},w_{i_{s+1}}}) =
(0,t/z_{f_e})$, $(z_{i_{n-1}},w_{i_{n-1}})=(1,-t/(1-z_{f_e}))$ and
$(z_n,w_n)=(\infty,0)$.
 Similarly,
for $|v_e|=3$, $(z_{i_{1}},w_{i_{1}})=(-t,1)$,
$(z_{i_{n-2}},w_{i_{n-2}}) = (0,\infty)$,
$(z_{i_{n-1}},w_{i_{n-1}})=(1,-t)$ and $(z_n,w_n)=(\infty,0)$. The
other labeled points are taken in an arbitrary position. The
component $z$ of the special points and the parameter $t$ provide a
coordinate system in $V \times B_\epsilon$.

Due to Knudsen's theorem there exists a unique $\kappa: V \times
B_\epsilon \to \cmod{n}$ which gives the family of $n$-pointed
curves given above.

\begin{lem}
\label{lem_projection1} %
$\mathrm{det} (d \kappa) \ne 0$ at $(\GS^o,\ve{p}^o) \in  D_\Gg$.
Hence, $\kappa$ gives  a local isomorphism.
\end{lem}

\begin{proof}
The parameter $t$ gives a regular function on $\kappa(V \times
B_\epsilon)$ which is vanishing along $D_\Gg \bigcap \kappa(V \times
B_\epsilon)$. The differential $d \kappa(\vec{v}) = \vec{v}$ for
$\vec{v} \in T_{(\GS^o,\ve{p}^o)}V$ since the restriction of
$\kappa$ on $V \times \{0\}$ is the identity map. We need to prove
that $d \kappa(\dd_t) \ne 0$. In other words, the curves are
non-isomorphic for different values of the parameter $t$. Let
$(\GS(t_i),\ve{p}(t_i)) \in V \times B_\epsilon$ be two $n$-pointed
curves for $t_1 \ne t_2$. A bi-holomorphic map $\GP: \GS(t_1) \to
\GS(t_2)$ is determined by the images of $p_{s+1},p_{i_{n-1}},p_n$
when $|v_e|>3$, and by the images of $p_{i_{n-2}},p_{i_{n-1}},p_n$
when $|v_e|=3$. However, the bi-holomorphic map $\GP$ mapping
$(p_{s+1},p_{i_{n-1}},p_n)(t_1) \mapsto
(p_{s+1},p_{i_{n-1}},p_n)(t_2)$ (resp. $(p_{i_{n-2}},
p_{i_{n-1}},p_n)(t_1) \mapsto (p_{i_{n-2}},p_{i_{n-1}},p_n)(t_2)$)
maps $p_{i_1}(t_1)=(z_{f_e}-t_1,1)$ to $(z_{f_e}-t_1,t_2/t_1) \ne
p_{i_1}(t_2)$ (resp.  $p_{i_1}(t_1)=(-t_1,1)$ to $(-t_1,t_2/t_1) \ne
p_{i_1}(t_2)$), i.e., $\GP$ can not be an isomorphism.
\end{proof}

\begin{rem} \label{rem_proje}
Due to Lemma \ref{lem_projection1}, the coordinates on $V \times
B_\epsilon$ provide a coordinate system at $(\GS^o;\ve{p}^o) \in
D_\Gg$. There is a natural coordinate projection  $\rho: V \times
B_\epsilon \to V$.

For a $\Gs$-invariant $\Gg$ and $c_\Gs$-invariant $V$, the above
coordinates and the local isomorphism $\kappa$ are equivariant with
respect to a suitable real structure ($(z,w) \mapsto
(\bar{z},\bar{w})$ when $\real{\GS} \ne \emp$, and $(z,w) \mapsto
(\bar{w},\bar{z})$ when $\real{\GS}= \emp$) on $\projc  \times
\projc$. Therefore, the real part $\real{V} \times ]-\epsilon,
\epsilon[$ of $V \times B_\epsilon$ provides a neighborhood for a
$(\GS^o;\ve{p}^o)$ in $\real{D_\Gg}$ with a set of coordinates on
it.
\end{rem}

\subsubsection{Contraction morphisms for o-planar trees.}
\label{sec_contraction_o} %

Let $(\Gg,\hat{o})$ be an o-planar tree and $\Gp: \Gg \to \Gt$ be a
morphism of $n$-trees contracting an invariant set of edges
$E_{con}= E_{\Gg} \smin \Gp_{E}(E_{\Gt})$. In such a situation, we
associate a particular o-planar structure $o$ on $\Gt$, as described
below in separate cases {\bf (a)} and {\bf (b)}, and speak of a {\it
contraction morphism} $\varphi: (\Gg,\hat{o}) \to (\Gt,o)$. In all
the cases, except {\bf (a-2)}, the o-planar structure $o$ is
uniquely defined by $\hat o$.

\begin{itemize}
\item[\bf (a)] Let $E_{con} = \{e=(f_e,f^e)\}$ and $e$ be an
invariant edge.

\begin{enumerate}
\item If $\dd_{\Gg}(e) =\{v_e,v^e\}\subset V_{\Gg}^\R$,
then we convert the o-planar structures
\begin{eqnarray*}
\begin{array}{lll}
\hat{o}_{v_e} &=& \{\R\GS_{v_e} \ne \emp; F^\pm_{\Gg}(v_e) ;
F_{\Gg}^\R (v_e) =\{ \{i_1\} <
\cdots <\{i_m\}<\{f_e\} \}\} \\
\hat{o}_{v^e} &=& \{\R\GS_{v^e} \ne \emp; F^\pm_{\Gg}(v^e) ;
F_{\Gg}^\R (v^e) =\{ \{i'_1\} < \cdots <\{i'_{m'}\}<\{f^e\} \}\}
\end{array}
\end{eqnarray*}
at $v_e$ and $v^e$  to an o-planar structure at vertex $v =
\Gp_{V}(\{v_e,v^e\})$ of $(\Gt,o)$ defining it by
\begin{eqnarray*}
\begin{array}{lll}
o_{v} &=& \{\R\GS_{v} \ne \emp; F^\pm_{\Gt}(v)
= F^\pm_{\Gg}(v_e)  \bigcup F^\pm_{\Gg}(v^e); \\
&& F_{\Gt}^\R (v) = \{ \{i_1\} < \cdots <\{i_m\}<\{i'_1\} < \cdots
<\{i'_{m'}\} \} \}.
\end{array}
\end{eqnarray*}
The o-planar structures are kept unchanged at all other invariant
vertices.

\item If $e$ is a special invariant edge, then we convert the o-planar
structure $\hat{o}=\{f_e \mapsto +, f^e \mapsto - \}$ of $\Gg$ into
an o-planar structure at the vertex $v= \Gp_{V}(\{v_e,v^e\})$ of
$\Gt$ defining it by
\begin{equation*}
o_v = \{\R\GS_{v} \ne \emp; F^+_{\Gt}(v) = F^+_{\Gg}(v_e) \smin
\{f_e\}, F^-_{\Gt}(v) = F^-_{\Gg}(v^e) \smin \{f^e\} ; F_{\Gt}^\R
(v) =\emp \}
\end{equation*}
or by
\begin{equation*}
o_v= \{\R\GS_{v} = \emp\}.
\end{equation*}
\end{enumerate}

\item[\bf (b)]  Let $E_{con}= \{e_i=(f_{e_i},f^{e_i}) \mid i=1,2\}$
where $f_{e_i}, i=1,2$ and $f^{e_i}, i=1,2$ are conjugate pairs of
flags.

\begin{enumerate}
\item If $\dd_{\Gg}(e_i) =\{\hat{v},v^{e_i}\}$, and $\hat{v} \in V_{\Gg}^\R$,
$v^{e_i} \not\in V_{\Gg}^\R$, then we convert the o-planar structure
\begin{eqnarray*}
o_{\hat{v}} = \{\R\GS_{\hat{v}} \ne \emp; F^\pm_{\Gg}(\hat{v}) ;
F_{\Gg}^\R (\hat{v}) =\{ \{i_1\} < \cdots <\{i_m\} \} \}.
\end{eqnarray*}
at $\hat{v}$ to an o-planar structure at
 $v=\Gp(\{\hat{v},v^{e_1},v^{e_2}\})$ of $\Gt$ defining it by
\begin{eqnarray*}
\begin{array}{lll}
o_v &=& \{\R\GS_{v} \ne \emp;
 F^+_{\Gt}(v) = F^+_{\Gg}(\hat{v})  \bigcup
 F_{\Gg}^+(v^{e_1}) \smin \{f_{e_1},f^{e_1} \}, \\
&&F^-_{\Gt}(v) = F^-_{\Gg}(\hat{v})  \bigcup
 F_{\Gg}^-(v^{e_2}) \smin \{f_{e_2},f^{e_2}\}; \\
&& F^\R_{\Gt}(v) =  \{ \{i_1\} < \cdots <\{i_m\} \} \}.
\end{array}
\end{eqnarray*}

\item If $\dd_{\Gg}(e_i) =\{\hat{v},v^{e_i}\}$, and $\hat{v} \in V_{\Gg}^\R$,
$v^{e_i} \not\in V_{\Gg}^\R$, then we convert the o-planar structure
$o_{\hat{v}} = \{\R\GS_{\hat{v}} = \emp \}$
at the vertex $\hat{v}$ to an o-planar structure at
 $v=\Gp(\{\hat{v},v^{e_1},v^{e_2}\})$ of $\Gt$ defining it by
$o_{v} = \{\R\GS_{v} = \emp \}.$

\item If $E_{con} =\{e_i=(f_{e_i},f^{e_i})\mid i=1,2\}$ and
 $\dd_{\Gg}(e_i) \bigcap V_{\Gg}^\R =\emp$, then we define
the o-planar structure at each $v$ in $\Gt$ to be the same as the
o-planar structure at $v$ of $(\Gg,\hat{o})$.
\end{enumerate}

\end{itemize}

\begin{rem}
Let $(\Gg,\hat{o})$ be an o-planar tree with $\real{\GS} \ne \emp$,
and let $\varphi_e: (\Gg,\hat{o}) \to (\Gt,o)$ be the contraction of
an edge $e \in E_{\Gg}$. If the o-planar tree $(\Gt,o)$ and the
u-planar tree $(\Gg,\hat{u})$ underlying $(\Gg,\hat{o})$ are given,
then the o-planar structure $\hat{o}$ can be reconstructed. For this
reason, when a stratum $\csq{\Gg,\hat{u}}$ contained in the boundary
of $\csqc{\Gt,u}$ is given, we denote the corresponding o-planar
structure $\hat{o}$ by $\Gd(o)$.
\end{rem}

\begin{prop}
\label{prop_boundary}%
A stratum $\csq{\Gg,\hat{u}}$ is contained in the boundary of
$\csqc{\Gt,u}$ if and only if the u-planar structures $u,\hat{u}$
can be lifted to o-planar structures $o,\hat{o}$ in such a way that
$(\Gt,o)$ is obtained by contracting an invariant set of edges of
$(\Gg,\hat{o})$.
\end{prop}

\begin{proof}
We need to consider  the statement only for the strata of
codimension one and two. These cases correspond to the contraction
morphisms from two/three-vertex o-planar (sub)trees to one-vertex
o-planar (sub)trees given in {\bf (a)} and {\bf (b)}. For a stratum
of higher codimension, the statement can be proved by applying the
elementary contractions {\bf (a)} and {\bf (b)} inductively. Here,
we consider only the case {\bf (a-1)}. The proof for other cases is
the same.

We first assume  that $(\Gt,o)$ is obtained by contracting the edge
$e$ of $(\Gg,\hat{o})$, where $(\Gg,\hat{o})$ is an o-planar
two-vertex tree with $V_\Gg = V_\Gg^\R = \{v_e,v^e\}$. An element
$\curve \in \csq{\Gg,\hat{o}}$ can be represented by the nodal curve
$\{(z-z_{f_e}) \cdot w =0\}$ in $\projc \times \projc$ with special
points $a_f = (z_f,0)$ and $b_f = (z_{f_e},w_f)$ such that
\begin{equation*}
\begin{array}{ll}
a_{f} \in \{ w =0 \ \& \ \Im(z)>0\}
& \mathrm{for}\  f \in F^+_{\Gg}(v_e)\\
a_{\bar{f}}  \in \{ w =0 \ \& \ \Im(z)<0\}
& \mathrm{for}\  \bar{f} \in F^-_{\Gg}(v_e) \\
\{ a_{i_1} < \cdots < a_{i_m}  \} \subset \{ w =0 \ \& \ \Im(z)=0\}
&  \mathrm{for}\  i_* \in F^\R_{\Gg}(v_e)
\end{array}
\end{equation*}
on the axis $w=0$, and $b_* = (z_{f_e},w_*)$
\begin{equation*}
\begin{array}{ll}
b_{f} \in \{ z -z_{f_e}=0 \ \& \ \Im(w)>0\}&
 \mathrm{for}\  f \in F^+_{\Gg}(v^e)\\
b_{\bar{f}} \in \{ z -z_{f_e}=0 \ \& \ \Im(w)<0\}&
\mathrm{for}\  \bar{f} \in F^-_{\Gg}(v^e) \\
\{ b_{i'_{1}} < \cdots < b_{i'_{m'}} \} \subset \{ z -z_{f_e}=0 \ \&
\ \Im(w)=0\}& \mathrm{for}\  i'_* \in F^\R_{\Gg}(v^e)
\end{array}
\end{equation*}
When we include the curve $\{(z-z_{f_e}) \cdot w = 0\}$ into the
family $\{(z-z_{f_e}) \cdot w +t =0\}$, the complex orientation
defined on the irreducible components $w=0$ and $z-z_{f_e}=0$ by the
halves $\Im(z)>0$ and, respectively, $\Im(w)>0$ extends continuously
to a complex orientation of $\{(z-z_{f_e}) \cdot w +t =0\}$ with
$t\in [0,\epsilon[$ defined by, say, $\Im(z)>0$. As a result, the
curves $\{(z-z_{f_e}) \cdot w +t =0\}$ with $t\in [0,\epsilon[$
acquire an o-planar structure given by
\begin{equation*}
\begin{array}{ll}
(z_f,w_f) \in \{ (z-z_{f_e}) \cdot w +t =0 \ \& \ \Im(z)>0\}
& \mathrm{for}\  f \in F^+_{\Gg}(v_e) \bigcup F^+_{\Gg}(v^e)\\
(z_{\konj{f}},w_{\konj{f}}) \in \{ (z-z_{f_e}) \cdot w +t=0 \ \& \
\Im(z)<0\}
& \mathrm{for}\ \bar{f} \in F^-_{\Gg}(v_e) \bigcup \in F^-_{\Gg}(v^e) \\
(z_f,w_f) \in \{(z-z_{f_e}) \cdot w +t=0 \ \& \ \Im(z)=0 \} &
\mathrm{for}\  f \in F^\R_{\Gg}(v_e) \bigcup F^\R_{\Gg}(v^e)
\end{array}
\end{equation*}
where the points on the real part of the curves $\{(z-z_{f_e})\cdot
w +t =0\}$ are cyclicly ordered by
\begin{equation*}
z_{i_1} <  \cdots < z_{i_m} <z_{i'_1} <  \cdots < z_{i'_{m'}}.
\end{equation*}
This is exactly the o-planar structure $(\Gt,o)$ defined in {\bf
(a-1)} of \ref{sec_contraction_o}.

Now assume that $\csq{\Gg,\hat{u}}$, where $(\Gg,\hat{u})$ is an
u-planar tree with $V_\Gg = V_\Gg^\R = \{v_e,v^e\}$, is contained in
the boundary of $\csqc{\Gt,u}$. There are four different o-planar
representatives of $(\Gg,\hat{u})$, and any two non reverse to each
other representatives $\hat o_1, \hat o_2$ provide by contraction
two different o-planar structures $(\Gt,o_i), i=1,2$. By the already
proved part of the statement, $\csq{\Gg,\hat{u}}$ is contained in
the boundary of $\csqc{\Gt,o_i}$ for each $ i=1,2$. It remains to
notice that any codimension one stratum is adjacent to at most two
main strata.
\end{proof}

\subsubsection{Examples}
\label{exa_realmoduli}

(i) The first nontrivial example is $\cmod{4}$. There are three real
structures: $c_{\Gs_1}$, $c_{\Gs_2}$, $c_{\Gs_3}$, where
\begin{eqnarray*}
\Gs_1 = \left(
\begin{array}{cccc}
1 & 2 & 3 & 4 \\
1 & 2 & 3 & 4 \\
\end{array}
\right), \ \Gs_2 = \left(
\begin{array}{cccc}
1 & 2 & 3 & 4 \\
2 & 1 & 3 & 4 \\
\end{array}
\right) \ \mathrm{and} \ \Gs_3 = \left(
\begin{array}{cccc}
1 & 2 & 3 & 4 \\
3 & 4 & 1 & 2 \\
\end{array}
\right).
\end{eqnarray*}
These real structures then give $\rmod{2k,l}$, where $(2k,l) =
(0,4), (2,2)$, and $(4,0)$ respectively.

\begin{figure}[htb]
\centerfig{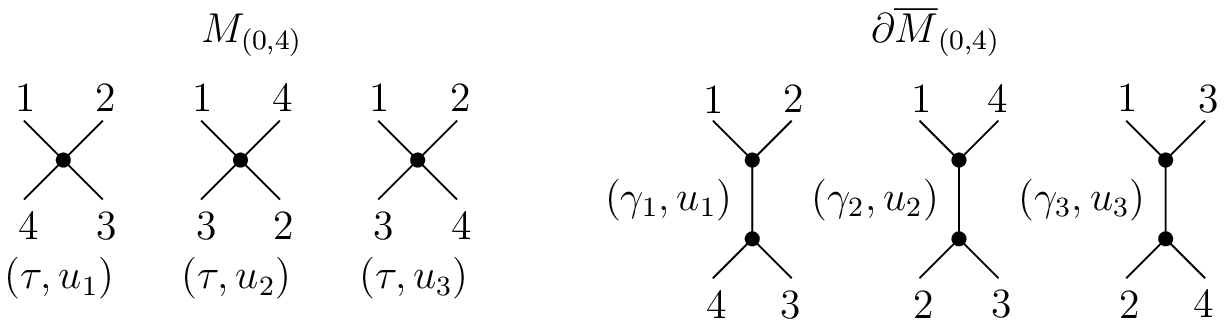,height=1.1in} %
\caption{U-planar trees encoding the strata of $\cmod{(0,4)}$.} %
\label{fig_rmoduli1}
\end{figure}

In the case $\Gs = id$, $\rmodo{0,4}$ is the configuration space of
four distinct points on $\projr$ up to the action of $PSL_{2}(\R)$.
The $4$-pointed curves $\curve \in \rmodo{0,4}$ can be identified
with $(0,x_2,1,\infty)$ where $x_2 \in \projr \smin \{0,1,\infty\}$.
Hence, $\rmodo{0,4} = \projr \smin \{0,1,\infty \}$ and its
compactification is $\rmod{0,4} = \projr$.  The three intervals of
$\rmodo{0,4}$ are the three configuration spaces $\csq{\Gt,u_i}$ and
the three  points  are the configuration spaces  $\csq{\Gg_i,u_i}$.
The u-planar trees $(\Gt,u_i)$ and $(\Gg_i,u_i)$ are given in Fig.
\ref{fig_rmoduli1}.

In the case $\Gs=\Gs_1$, $\rmodo{2,2}$ is the space of distinct
configurations of two points in $\projr$ and a pair of complex
conjugate points in $\projc \smin \projr$. $\curve \in \rmodo{2,2}$
is identified with $(\im,-\im, x_3, \infty) \in \csq{\Gt,u},
 -\infty < x_3 <\infty$. Hence, $\rmodo{2,2}= \projr \smin \{\infty\}$
and its compactification is $\rmod{2,2} = \projr$. The interval
$\rmodo{2,2}$ is $\csq{\Gt,u}$ and the  point at its closure is
$\csq{\Gg,u}$.

In the case $\Gs=\Gs_2$, the space $\rmodo{4,0}$ has different
pieces parameterizing real curves with non-empty and empty  real
parts: The subspace of $\rmodo{4,0}$ parameterizing the
$(4,0)$-pointed real curve with $\real{\GS} \not= \emp$ is  $(\Gl
\im, \im, -\Gl \im, -\im)$ where $\Gl \in ]-1,1[ \smin \{0\}$. The
subspace of $\rmodo{4,0}$ parameterizing the real curves with
$\real{\GS}= \emp$ is $(\Gl \im, \im, - \im/ \Gl,- \im)$, where $\Gl
\in ]-1,1[$. Note that, the pieces parameterizing $\real{\GS} \not=
\emp$ and $\real{\GS} = \emp$ are joined through the boundary points
corresponding to curves with isolated real singular points. The
compactification $\rmod{4,0}$ is $\projr$.

(ii) The moduli space $\cmod{5}$  has three different real
structures $c_{\Gs_1}, c_{\Gs_2}$ and $c_{\Gs_3}$ where
\begin{eqnarray*}
\Gs_1 = \left(
\begin{array}{ccccc}
1 & 2 & 3 & 4 &5 \\
1 & 2 & 3 & 4 &5 \\
\end{array}
\right), \ \Gs_2 = \left(
\begin{array}{ccccc}
1 & 2 & 3 & 4 & 5\\
2 & 1 & 3 & 4 & 5\\
\end{array}
\right) \ \mathrm{and} \ \Gs_3 = \left(
\begin{array}{ccccc}
1 & 2 & 3 & 4 & 5 \\
3 & 4 & 1 & 2 & 5\\
\end{array}
\right).
\end{eqnarray*}

\begin{figure}[htb]
\centerfig{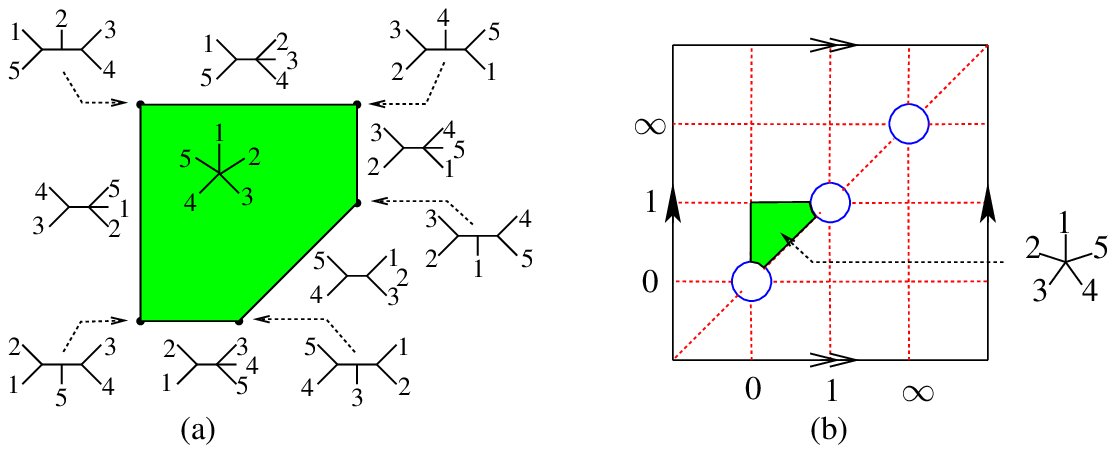,height=2 in} %
\caption{(a) Stratification of  $\cspq{\Gt,u}$.
(b) The stratification of $\rmod{0,5}$.} %
\label{fig_rmoduli3}
\end{figure}

The space $\rmodo{0,5}$ is identified with the configuration space
of five distinct points on $\projr$ modulo $PSL_{2}(\R)$. It is
$(\projr \setminus \{0,1,\infty\})^{2} \setminus \GD$, where $\GD$
is union of all diagonals. Each connected component of $\rmodo{0,5}$
is isomorphic to a two dimensional simplex.  The closure of each
cell can be obtained by adding the boundaries given in Section
\ref{sec_stratification}; for an example see Fig.
\ref{fig_rmoduli3}a. It  gives the compactification of $\rmod{0,5}$
which is  a torus with 3 points blown up: the cells corresponding to
u-planar trees $(\Gt,u_1)$ and $(\Gt,u_2)$ are glued along the face
corresponding to $(\Gg,u)$ which gives $(\Gt,u_i),i=1,2$ by
contracting some edges, see Fig. \ref{fig_rmoduli3}b.

The space $\rmodo{2,3}$ is isomorphic to configurations of a
conjugate pair of points on $\projc$. The automorphisms allows us to
identify such configurations with $(z,\bar{z},0,1,\infty)$ where $z
\in \C \smin \R$. Hence, it can be given as $\projc \smin \projr$.
The $\rmod{2,3}$  is obtained as a sphere with 3 points blown up
according to the stratification given in  Section
\ref{sec_stratification}.

Finally, elements of $\rmodo{4,1}$ can be identified with
$(z,\im,\bar{z},-\im,\infty)$. Hence it can be identified with
$\projc \smin (\projr \cup \{\im,-\im\})$. Therefore, connected
components are isomorphic to $\H^+ \smin \{\im\}$. The $\rmod{4,1}$
is a sphere with a point blown up.

The moduli space $\cmod{5}$ is a del Pezzo surface of degree 5 and
these are all the possible real parts of this del Pezzo surface (see
\cite{dik}).

\section{The first Stiefel-Whitney class of $\rmod{2k,l}$}
\label{sec_steifel} %

In this section we calculate the first Stiefel-Whitney class of
$\rmod{2k,l}$ by using the stratification given in  Theorem
\ref{thm_real_strata}.

\subsection{Orientations of top-dimensional strata}
\label{sec_ori}%
Let $(\Gt,o)$ be  a one-vertex o-planar tree. The coordinates on
$\csq{\Gt,o}$ given in  Section  \ref{sec_normal_pos} determine an
orientation of $\csq{\Gt,o}$. For instance, let $|\fix(\Gs)| \geq 3$
and let the o-planar structure on $(\Gt,o)$  be given by $\Conj^\pm$
and by a linear ordering $x_{i_1}=0 < x_{i_2} <\cdots < x_{i_{l-1}}
=1 < x_{i_l} := x_n =\infty$ on $\fix(\Gs)$. The  coordinates  in
{\bf (A)} of \ref{sec_normal_pos} generate the following
top-dimensional differential form on $\csq{\Gt,o}$:
\begin{eqnarray} \label{eqn_orient}
\form{\Gt,o} := \left(\frac{\im}{2} \right)^k \bigwedge_{\Ga_* \in
\Conj^+ } dz_{\Ga_*} \wedge d\overline{z}_{\Ga_*} \ \bigwedge \
dx_{i_2} \wedge \cdots \wedge dx_{i_{l-2}}.
\end{eqnarray}
The multiplication of top-dimensional forms with a positive valued
function $\Theta: \csq{\Gt,o} \to \R_{>0}$ defines an equivalence
relation on sections of $\det(T \csq{\Gt,o})$. An {\it orientation}
is an equivalence class of nowhere zero top-dimensional forms with
respect to this equivalence relation. We denote the equivalence
class of $\form{\Gt,o}$ by  $\ori{\Gt,o}$.

Similarly, using the  coordinates given in  {\bf (A$'$)}, {\bf (B)},
{\bf (C)} and {\bf (D)} in \ref{sec_normal_pos} and their ordering,
we determine  differential forms $\Go_{(\Gt,o)}$ and  orientations
$\ori{\Gt,o}$ of  $\csq{\Gt,o}$ for all $(\Gt,o)$ with
$|V_{\Gt}|=1$.

\subsection{Orientations of codimension one strata}
\label{sec_d_ori} Let $(\Gg,o)$ be a two-vertex o-planar tree. Let
$V_\Gg =\{v_e,v^e\}$ and $e = (f_e, f^e)$ be the edge where
$\dd_\Gg(n) = \dd_\Gg(f_e) = v_e$ and $\dd_\Gg(f^e) = v^{e}$.

By choosing three flags in $F_\Gg(v_e)$ and $F_\Gg(v^e)$, and using
the calibrations as in \ref{sec_normal_pos} we obtain a coordinate
system in $\csq{\Gg_v,o_v}$ for each $v \in \{v_e,v^e\}$.  More
precisely, we use the following choice.

\begin{itemize}
\item[\bf I.] Let $\real{\GS_v} \ne \emp$ for
$v \in \{v_e,v^e\}$ and let $\fix(\Gs) \ne \emp$. If
$|F_\Gg^\R(v_e)| \geq 3$ (resp. $|F_\Gg^\R(v^e)| \geq 3$), then we
specify an isomorphism $\Phi_{v_e}: \GS_{v_e} \to \projc$ (resp.
$\Phi_{v^e}: \GS_{v^e} \to \projc$) by mapping three consecutive
special points as follows: If $|F_\Gg^\R(v_e)| > 3$ and the special
points $p_{f_e}$ and $p_n$ are not consecutive, then $\GP_{v_e}:
(p_{i_{l-1}},p_{n}, p_{i_{1}}) \mapsto (1,\infty,0)$. If
$|F_\Gg^\R(v_e)| \geq 3$ and the special points $p_{f_e}$ and $p_n$
are consecutive and
\begin{eqnarray*}
\begin{array}{ll}
\{f_e\} < \{n\} < \{i_1\}, & \Longrightarrow   \\
\{i_{l-1}\} < \{n\} < \{f_e\}, & \Longrightarrow\
\end{array}
\begin{array}{lcl}
\GP_{v_e}: (p_{f_e},p_{n}, p_{i_{1}}) &\mapsto& (1,\infty,0), \\
\GP_{v_e}: (p_{i_{l-1}},p_{n}, p_{f_e}) &\mapsto& (1,\infty,0).
\end{array}
\end{eqnarray*}
For $|F_\Gg^\R(v^e)| \geq 3$, $\GP_{v^e}:
(p_{i_{q+1}},p_{i_{q+r}},p_{f^e}) \mapsto (1,\infty,0)$.

If $|F_\Gg^\R(v_e)| < 3$ (resp. $|F_\Gg^\R(v^e)| < 3$), then in
addition to $p_n \mapsto \infty$ (resp. $p_{f^{e}} \mapsto 0$), we
pick the maximal element $\Ga$ in $F^+_\Gg(v_e)$ (resp.
$F^+_\Gg(v^e)$) and map the pair of  conjugate labeled points
$(p_\Ga,p_{\konj{\Ga}})$ to $(\im,-\im)$.

\item[\bf II.] Let $\real{\GS_v} \ne \emp$ for
$v \in \{v_e,v^e\}$ and let $\fix(\Gs) =\emp$. We specify an
isomorphism  $\Phi_v :\GS_{v} \to \projc$ by mapping the pair of
conjugate labeled points $(p_\Ga,p_{\konj{\Ga}})$  to $(\im,-\im)$
for the maximal element $\Ga$ in $F_\Gg^+(v)$, and $p_{f_e} \mapsto
0$ (resp. $p_{f^e} \mapsto 0$).

\item[\bf III.]  Let $\real{\GS}$ be an isolated
real node. We pick a maximal element $\Ga_{k-1}$ in $F_\Gg^+(v_e)
\smin \{n\}$ and specify isomorphisms $\Phi_{v_e} :\GS_{v_e} \to
\projc$ and $\Phi_{v^e} :\GS_{v^e} \to \projc$ by mapping the
special points $(p_{f_e},p_{\Ga_{k-1}},p_{n})$ and, respectively,
$(p_{f^e},p_{\konj{\Ga}_{k-1}},p_{\konj{n}})$ to $(0,\im/2,\im)$.
\end{itemize}

By using the o-planar structure
\begin{equation*}
o_v= \left\lbrace
\begin{array}{ll}
\{\real{\GS_v} \ne \emp; F_\Gg^\pm(v); F_\Gg^\R(v) =
\{\{f_1\}<\cdots<\{f_{l_v}\} \} \} &
\mathrm{for\ case} \ \mathbf{I}, \\
\{\real{\GS_v} \ne \emp; F_\Gg^\pm(v); F_\Gg^\R(v) = \emp \} &
\mathrm{for\ case} \ \mathbf{II},
\end{array}
\right.
\end{equation*}
of $\Gg_v$ for each $v \in \{v_e,v^e\}$, arrange the coordinates of
the special points in the following order
\begin{equation*}
(z_{\Ga_1}, \cdots, z_{\Ga_{k_{v}}}, x_{i_1},\cdots,x_{i_{l_v}}),
\end{equation*}
where $\Ga_* \in F^+_\Gg(v)$. We fix special points as in ({\bf I})
and ({\bf II}), and apply (\ref{eqn_orient}) to introduce
top-dimensional differential forms $\GO_{(\Gg_{v^e},o_{v^e})}$ and
$\GO_{(\Gg_{v_e},o_{v_e})}$ on $\csq{\Gg_{v^e},o_{v^e}}$ and
$\csq{\Gg_{v_e},o_{v_e}}$ (note that the resulting forms do not
depend on the order of $z$-coordinates). In the case ({\bf III}),
there are no real special points, so we may get a top-dimensional
differential form
$\GO_{(\Gg,o)}$ on $\csq{\Gg,o}$ 
via choosing the vertex $v \in V^+_\Gg$  with ordering arbitrarily
the $z$-coordinates
\begin{equation*}
(z_{\Ga_1}, \cdots, z_{\Ga_{k_{v}}})
\end{equation*}
where $F_\Gg^+ =\{\Ga_1,\cdots,\Ga_{k_v} \}$.

In such a way, we produce well-defined orientations
$\orie{\Gg_{v^e},o_{v^e}}$ and  $\orie{\Gg_{v_e},o_{v_e}}$ of,
$\csq{\Gg_{v^e},o_{v^e}}$ and $\csq{\Gg_{v_e},o_{v_e}}$ respectively, 
and finally get an orientation on $\csq{\Gg,o}$ determined by
\begin{eqnarray*} \label{eqn_bd_or}
\begin{array}{ll}
\orie{\Gg_{v^e},o_{v^e}} \wedge \orie{\Gg_{v_e},o_{v_e}}
& \mathrm{when}\ V_\Gg=V_\Gg^\R                           \\
\orie{\Gg,o} & \mathrm{when}\ V_\Gg^\R=\emp \ \mathrm{and}\ v \in
V_\Gg^+.
\end{array}
\end{eqnarray*}

\subsection{Induced orientations on codimension one strata}
\label{sec_in_ori}%
Let $\csq{\Gt,u}$ be a top dimensional stratum, and let
$\csq{\Gg,\hat{u}}$ be a codimension one stratum contained in the
boundary of $\csqc{\Gt,u}$. We lift the u-planar structures
$u,\hat{u}$ to o-planar representatives $o, \hat{o}=\Gd(o)$ such
that $(\Gt,o)$ is obtained by contracting the edge of $(\Gg,\Gd(o))$
(see Prop. \ref{prop_boundary}). Then we pick a point
$(\GS^o,\ve{p}^o) \in \csq{\Gg,\Gd(o)}$ and consider a tubular
neighborhood $\real{V} \times [0,\epsilon[$ of $(\GS^o,\ve{p}^o)$ in
$\csqc{\Gt,o}$ as in Section \ref{sec_stratification}.

The orientation $\ori{\Gt,o}$, introduced in \ref{sec_ori}, induces
some orientation on $\csq{\Gg,\Gd(o)}$: The outward normal direction
of $\csqc{\Gt,o}$ on $\real{V} \times \{0\} \subset
\csq{\Gg,\Gd(o)}$ is $-\dd_t$, where $t$ is the standard coordinate
on $[0,\epsilon[\subset \R$. Therefore a differential form
$\form{\Gg,\Gd(o)}$ defines the induced orientation, if and only if
\begin{equation} \label{eqn_form}
- dt \wedge \Go_{(\Gg,\Gd(o))} = \Theta \, \form{\Gt,o}
\end{equation}
with $\Theta >0$ at each point of $\real{V} \times ]0,\epsilon[$.

In what follows, we compare the induced orientation
$\ori{\Gg,\Gd(o)}$ with $\orie{\Gg_{v^e},o_{v^e}} \wedge
\orie{\Gg_{v_e},o_{v_e}}$.

\subsubsection{Case I: $\fix(\Gs) \ne \emp$.}

\begin{lem} \label{lem_detay0}
Let $(\Gg,\Gd(o))$ be an o-planar tree as above, and let
$|F^\R_\Gg|=l+2>3$, and $|F^\R_\Gg(v^e)|=r+1$. Then,
\begin{eqnarray}
\ori{\Gg,\Gd(o)} = (-1)^{\GH}\ \orie{\Gg_{v^e},o_{v^e}} \wedge
\orie{\Gg_{v_e},o_{v_e}} \nonumber
\end{eqnarray}
where the values of $\aleph$ for separate cases are given in the
following table.

\begin{tabular}{|c|c|c|c|c|}
\hline
$\GH$ & $l-r \geq 3$ & \multicolumn{2}{|c|}{$l-r =2$} &  $l-r =1$  \\
\cline{3-4}
       &              & $\{i_{1}\}<\{f_e\}<\{n\}$&$\{f_e\}<\{i_{l-1}\}< \{n\}$&    \\
\hline
$r \geq 2$ &$(q+1)(r+1)$ &         $0$           &         $l+1$        &   $l+1$    \\
\hline
$r   =  1$ &    $1$      &         $1$           &          $1$         &    $1$     \\
\hline
$r   =  0$ &   $q+1$     &         $0$           &          $0$         &    $0$     \\
\hline
\end{tabular}

\noindent Here, the third and  fourth columns correspond to two
possible cyclic orderings of $F_\Gg^\R(v_e)$ for $|F_\Gg^\R(v_e)|=3$
in  Case {\bf I} of Section \ref{sec_d_ori}.
\end{lem}

\begin{figure}[htb]
\centerfig{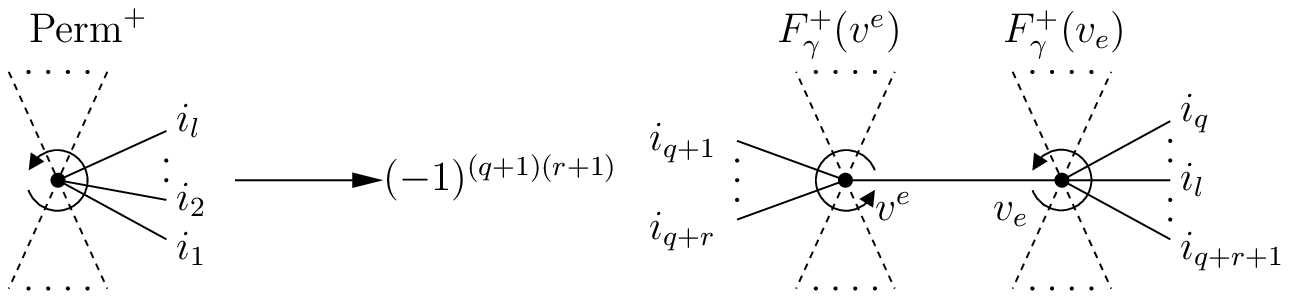,height=1in} %
\caption{Codimension 1 boundaries of $\csqc{\Gt,\Gd(o)}$ where $l-r
\geq 3$
\& $r \geq 2$.} %
\label{fig_orientation}
\end{figure}

\begin{proof}
We will prove the statement only for the special case of $l-r \geq
3$, $r \geq 2$.  The calculations for other cases are almost
identical.

Let $(\GS^o,\ve{p}^o) \in \csq{\Gg,\Gd(o)}$. We set $\GS_{v_e}^o$ to
be $\{w=0\}$ and $\GS_{v^e}^o$ to be $\{z -x_{f_e}=0\}$. According
to the convention in  \ref{sec_d_ori}, the consecutive special
points $(p_{i_{l-1}},p_{n},p_{i_{1}})$ (resp.
$(p_{i_{q+1}},p_{i_{q+r}},p_{f^e})$) on the component $\GS_{v_e}$
(resp. $\GS_{v^e}$) are fixed at $(1,\infty,0)$. As shown in the
proof of Proposition \ref{prop_boundary}, a tubular neighborhood
$\real{V} \times [0,\epsilon[$ of $(\GS^o,\ve{p}^o)$ in
$\csqc{\Gt,o}$ can be given by the family $\{(z-x_{f_e}) \cdot w +t
=0 \mid t\in [0, \epsilon[ \}$ with labeled  points $p_{i_{1}}=
(0,t/x_{f_e})$, $p_{i_{q+1}}= (x_{f_e}-t,1)$, $p_{i_{q+r}}=
(x_{f_e}, \infty)$, $p_{i_{l-1}}= (1,-t/(1 -x_{f_e}))$,
$p_{n}=(\infty,0)$,
 $p_{i_*}=(x_{i_*},-t/(x_{i_*} -x_{f_e}))$ for
$i_* \in F_\Gg^\R \smin \{i_1,i_{q+1},i_{q+r},i_{l-1},n\}$ and
$p_\Ga= (z_\Ga,-t/(z_{\Ga} -x_{f_e}))$ for $\Ga \in F_\Gg^+$.

We first consider the following subcase: the special points
$p_{f_e}$ and $p_n$ are not consecutive. According to
\ref{sec_d_ori}, the differential forms $\GO_{(\Gg_{v_e},o_{v_e})}$
and $\GO_{(\Gg_{v^e},o_{v^e})}$ of this case are as follows:
\begin{eqnarray*}
&\GO_{(\Gg_{v_e},o_{v_e})} &=  \left(
\frac{\im}{2}\right)^{|F_\Gg^+(v_e)|} \bigwedge_{\Ga \in
F^+_\Gg(v_e)}
d z_\Ga \wedge d \konj{z}_\Ga \wedge \\
&& dx_{i_2} \wedge \cdots \wedge dx_{i_q} \wedge dx_{f_e} \wedge
\widehat{dx_{i_{q+1}}} \wedge \widehat{\cdots} \wedge
\widehat{dx_{i_{q+r}}}
\wedge dx_{i_{q+r+1}} \wedge \cdots \wedge dx_{i_{l-2}}, \\
&\GO_{(\Gg_{v^e},o_{v^e})} &=  \left(
\frac{\im}{2}\right)^{|F_\Gg^+(v^e)|} \bigwedge_{\Gb \in
F^+_\Gg(v^e)} d w_\Gb \wedge d \konj{w}_\Gb \bigwedge dy_{i_{q+2}}
\wedge \cdots \wedge dy_{i_{q+r-1}}
\end{eqnarray*}

By using the  identities $w_\Gb = -t/(z_\Gb -x_{f_e})$ for $\Gb \in
F^+_\Gg(v^e)$ and $y_i = -t/(x_i - x_{f_e})$ for $i \in
F_\Gg^\R(v^e)$,  we obtain the following equalities:
\begin{eqnarray*}
\begin{array}{ll}
d t      = - d x_{i_{q+1}} + dx_{f_e}, \ \
dx_{f_e} =   d x_{i_{q+r}}, & \\
d w_\Gb =  -\frac{dt}{z_\Gb -x_{f_e}}+ \frac{t
dz_\Gb}{(z_\Gb-x_{f_e})^2}
- \frac{t dx_{f_e}}{(z_\Gb-x_{f_e})^2} & \mathrm{for}\ \Gb \in F_\Gg^\pm(v^e), \\
d y_i   =  -\frac{dt}{x_i -x_{f_e}}+ \frac{t dx_i}{(x_i-x_{f_e})^2}
- \frac{t dx_{f_e}}{(x_i-x_{f_e})^2}        & \mathrm{for}\
i  =q+2,\cdots,q+r-1. \\
\end{array}
\end{eqnarray*}
These identities imply that $- dt \wedge \GO_{(\Gg_{v^e},o_{v^e})}
\wedge \GO_{(\Gg_{v_e},o_{v_e})}$ is equal to
\begin{eqnarray*}
(-1)^{(r-1)(q-1)} \left( \frac{\im}{2}\right)^{|F_\Gg^+|}  \Theta
\bigwedge_{\Ga \in F_\Gg^+} d z_\Ga \wedge d \konj{z}_\Ga \bigwedge
dx_{i_2} \wedge  \cdots \wedge dx_{l-2}
\end{eqnarray*}
where $\Theta =  \prod_{\Gb \in F_\Gg^+(v^e)} t (z_\Gb
-x_{f_e})^{-2} \prod_{i_{q+2},\cdots,i_{q+r-1}} t (x_i
-x_{f_e})^{-2}$. Since $\Theta >0$, the orientation defined by $-dt
\wedge \GO_{(\Gg_{v^e},o_{v^e})} \wedge \GO_{(\Gg_{v^e},o_{v^e})}$
is equal to $(-1)^\GH \ori{\Gt,o}$.

We now consider the cases $\{f_e\}<\{n\}<\{i_{1}\}$  (i.e.,
$q+r=l-1$) and $\{i_{l-1}\}<\{n\}<\{f_e\}$ (i.e., $q=0$). According
to \ref{sec_d_ori}, the differential forms
$\GO_{(\Gg_{v^e},o_{v^e})} \wedge \GO_{(\Gg_{v_e},o_{v_e})} $ are
equal to
\begin{eqnarray*}
\left( \frac{\im}{2}\right)^{|F_\Gg^+|} & \left( \bigwedge_{\Gb \in
F^+_\Gg(v^e)} d z_\Gb \wedge d \konj{z}_\Gb \bigwedge
dy_{i_{q+2}} \wedge \cdots \wedge dy_{i_{l-2}}\right) \\
& \wedge \left( \bigwedge_{\Ga \in F^+_\Gg(v_e)} d z_\Ga \wedge d
\konj{z}_\Ga \bigwedge dx_{i_{2}} \wedge \cdots \wedge
dx_{i_q}\right)
\end{eqnarray*}
when $q+r=l-1$, and
\begin{eqnarray*}
\left( \frac{\im}{2}\right)^{|F_\Gg^+|} & \left( \bigwedge_{\Gb \in
F^+_\Gg(v^e)} d z_\Gb \wedge d \konj{z}_\Gb \bigwedge
dy_{i_2} \wedge \cdots \wedge dy_{i_{r-1}}\right) \\
& \wedge \left( \bigwedge_{\Ga \in F^+_\Gg(v_e)} d z_\Ga \wedge d
\konj{z}_\Ga \bigwedge dx_{i_{r+1}} \wedge \cdots \wedge
dx_{i_{l-2}}\right)
\end{eqnarray*}
when $q=0$. The equation $(z-x_{f_e}) \cdot w +t =0$ implies the following
equalities:
\begin{eqnarray*}
\begin{array}{ll}
d t      = - d x_{i_{q+1}}, \ \ dx_{f_e} =   d x_{i_{l-1}}
& \mathrm{when} \ q+r=l-1  \\
d t      =  d x_{i_r}, \ \ dx_{f_e} =   d x_{i_{q+r}}
& \mathrm{when} \ q=0,  \\
\end{array}
\end{eqnarray*}
and
\begin{eqnarray*}
\begin{array}{ll}
d w_\Gb =  -\frac{dt}{z_\Gb -x_{f_e}}+ \frac{t
dz_\Gb}{(z_\Gb-x_{f_e})^2}
- \frac{t dx_{f_e}}{(z_\Gb-x_{f_e})^2} & \mathrm{for}\ \Gb \in F_\Gg^\pm(v^e), \\
d y_i   =  -\frac{dt}{x_i -x_{f_e}}+ \frac{t dx_i}{(x_i-x_{f_e})^2}
- \frac{t dx_{f_e}}{(x_i-x_{f_e})^2}        & \mathrm{for}\
i  =i_{q+2},\cdots,i_{q+r-1}. \\
\end{array}
\end{eqnarray*}
By using these identities we obtain that $- dt \wedge
\GO_{(\Gg_{v^e},o_{v^e})} \wedge \GO_{(\Gg_{v^e},o_{v^e})}$ is equal
to
\begin{eqnarray*}
(-1)^{(q-1)(l-q-2)} \left( \frac{\im}{2}\right)^{|F_\Gg^+|}  \Theta
\bigwedge_{\Ga \in F_\Gg^+} d z_\Ga \wedge d \konj{z}_\Ga \bigwedge
dx_{i_2} \wedge  \cdots \wedge dx_{i_{l-2}}
\end{eqnarray*}
when $q+r=l-1$, and
\begin{eqnarray*}
(-1)^{(r-1)} \left( \frac{\im}{2}\right)^{|F_\Gg^+|}  \Theta
\bigwedge_{\Ga \in F_\Gg^+} d z_\Ga \wedge d \konj{z}_\Ga \bigwedge
dx_{i_2} \wedge  \cdots \wedge dx_{i_{l-2}}
\end{eqnarray*}
when $q=0$. Since $\Theta = \prod_{\Gb \in F_\Gg^+(v^e)} t (z_\Gb
-x_{f_e})^{-2} \prod_{i_{q+2},\cdots,i_{q+r-1}} t (x_i
-x_{f_e})^{-2}  >0$, the orientation $\ori{\Gg,\Gd(o)}$ induced by
$\ori{\Gt,o}$ is equal to
\begin{eqnarray*}
 &(-1)^\GH&  \left[ \GO_{(\Gg_{v^e},o_{v^e})} \wedge
\GO_{(\Gg_{v_e},o_{v_e})} \right]  = \\
&& \left\lbrace
\begin{array}{ll}
(-1)^{(q+1)(r-1)} \left[ \GO_{(\Gg_{v^e},o_{v^e})} \wedge
\GO_{(\Gg_{v_e},o_{v_e})}\right]
 & \mathrm{when} \ q+r=l-1,\\
(-1)^{(r-1)} \left[ \GO_{(\Gg_{v^e},o_{v^e})} \wedge
\GO_{(\Gg_{v_e},o_{v_e})}\right]
 & \mathrm{when} \ q=0.
\end{array}\right.
\end{eqnarray*}
\end{proof}

\subsubsection{Case II.\   $\fix(\Gs)=\emp$.}
The  different cases for boundaries of $\csq{\Gt,o}$  are treated
separately. The proofs are essentially the same as the proof of
Lemma \ref{lem_detay0}.

\paragraph{\bf Subcase  $\real{\GS} \ne \emp$.}
Let $(\Gt,o)$ be a one-vertex o-planar tree with $\real{\GS} \ne
\emp$, and let $\ori{\Gt,o}$ be, in accordance with \ref{sec_ori},
the orientation of $\csq{\Gt,o}$ defined by the differential form
\begin{equation} \label{eqn_gerek}
\form{\Gt,o} := (\frac{\im}{2})^{k-2} \bigwedge_{ \Ga_* \in  \Conj^+
\smin (\Conj^+ \bigcap \{k-1,k,2k-1,2k\})} dz_{\Ga_{*}} \wedge
d\overline{z}_{\Ga_{*}} \bigwedge d\Gl
\end{equation}
(which is given  by the coordinates in {\bf (C)} of
\ref{sec_normal_pos}). Here $\Gl = \Im(z_{\Ga'})$ and $\Ga' \in
\{k-1,2k-1\} \bigcap \Conj^+$.

\begin{lem} \label{lem_detay1}
Let $(\Gg,\Gd(o))$ be a two-vertex o-planar tree, and let the
corresponding strata $\csq{\Gg,\Gd(o)}$ be contained in the boundary
of $\csqc{\Gt,o}$.

(i) If $V_\Gg^\R = V_\Gg$, then the orientation $\ori{\Gg,\Gd(o)}$
induced by the orientation $\ori{\Gt,o}$ is equal to
$-\orie{\Gg_{v^e},o_{v^e}} \wedge \orie{\Gg_{v_e},o_{v_e}}.$

(ii) If $V_\Gg^\R = \emp$, then the orientation $\ori{\Gg,\Gd(o)}$
induced by the orientation $\ori{\Gt,o}$ is equal to
$\orie{\Gg,\Gd(o)}$.
\end{lem}

\paragraph{\bf Subcase  $\real{\GS} = \emp$.}
Let $(\Gt,o)$ be a one-vertex o-planar tree where
$o=\{\real{\GS}=\emp\}$ and let $\ori{\Gt,o}$ be, in accordance with
\ref{sec_ori}, the orientation of $\csq{\Gt,o}$ defined by
\begin{equation} \label{eqn_gerek2}
\form{\Gt,o} := -(\frac{\im}{2})^{k-2} \bigwedge_{\Ga_* \in
\{1,\cdots,k-2\}} dz_{\Ga_{*}} \wedge d\overline{z}_{\Ga_{*}}
\bigwedge d\Gl
\end{equation}
(which is given  by the coordinates in {\bf (D)} of
\ref{sec_normal_pos}). Here $\Gl = \Im(z_{k-1})$.

\begin{lem} \label{lem_detay2}
Let $(\Gg,\Gd(o))$ be a two-vertex o-planar tree where $V_\Gg^\R =
\emp$, and let $\csq{\Gg,\Gd(o)}$ be contained in the boundary  of
strata $\csqc{\Gt,o}$ given above. Then the orientation
$\ori{\Gg,\Gd(o)}$ induced by the orientation $\ori{\Gt,o}$ is equal
to $(-1)^\aleph \orie{\Gg,\Gd(o)}$ where $\aleph=|\{1,\cdots,k-1\}
\bigcap F_\Gg^-|+1$.
\end{lem}

\subsection{Conventions} \label{sec_convention}
Let $(\Gt,o_\star)$ be the one-vertex o-planar tree where the
o-planar structure $o_\star$ is given  by $\real{\GS} \ne \emp$,
$F_{\Gt}^+ =\{1,2,\cdots,k\}$, $F_{\Gt}^- =\{k+1,\cdots,2k\}$, and
$F_{\Gt}^\R = \{2k+1\} < \{2k+2\}< \cdots < \{2k+l:=n\}$. All the
other o-planar structures with $\real{\GS}\ne \emp$ on $\Gt$ are
obtained as follows.

Let $\Gr \in S_n$ be a permutation which commutes with $\Gs$ and, if
$l>0$, preserves $n$. It determines an o-planar structure given by
$\Gr(o_\star)= \{ \real{\GS}\ne\emp; \Gr(\Conj^\pm); \fix(\Gs)
=\{\{\Gr(2k+1)\}< \cdots <\{\Gr(2k+l-1)\}< \{\Gr(n)=n\}\}\}$. The
parity of $\Gr$ depends only on $o=\Gr(o_\star)$ and we call it {\it
parity}  $|o|$ of $o=\Gr(o_\star)$.

\subsubsection{Convention of orientations} \label{o_convention}
We fix an orientation for each top-dimensional stratum as follows.

\begin{itemize}
\item[{\bf a.}] {\bf Case $\real{\GS} \ne \emp$}.
First, we select o-planar representatives for each one-vertex
u-planar tree with $\real{\GS} \ne \emp$ as follows:

\begin{enumerate}
\item if $l \geq 3$, we choose the representative
$(\Gt,o)$ of $(\Gt,u)$ for which $\{2k+1\} < \{n-1\} < \{n\}$;

\item
if $l < 3$, we choose the representative $(\Gt,o)$ of $(\Gt,u)$ such
that $k \in \Conj^+$.
\end{enumerate}

\noindent We denote the set of these o-planar representatives of
u-planar trees by $\otree$, and select the orientation for
$\csq{\Gt,u}=\csq{\Gt,o}$ with $\csq{\Gt,o}\in\otree$ to be
\begin{equation} \label{eqn_or}
(-1)^{|o|} \ori{\Gt,o},
\end{equation}
where $\form{\Gt,o}$ is the form defined according to \ref{sec_ori}
and $|o|$ is the parity introduced in \ref{sec_convention}.

\item[{\bf b.}] {\bf Case $\real{\GS} = \emp$.}
Here, we choose the orientation defined by the form
(\ref{eqn_gerek2}).
\end{itemize}

In what follows, if $\real{\GS} \ne \emp$ we denote the set of flags
$\{2k+1,n-1,n\}$ (for $l \geq 3$ case) and $\{k,2k,n\}$ (for $l<3$
case) by $\fF$.

\subsection{Adjacent top-dimensional strata with $\real{\GS} \ne \emp$}
\label{sec_ad}%
Let $\csq{\Gt,u_i},i=1,2,$ be a pair of adjacent top-dimensional
strata with $\real{\GS} \ne \emp$, and $\csq{\Gg,u}$ be their common
codimension one boundary stratum. Let $(\Gt,o_i)$ be the o-planar
representatives of $(\Gt,u_i)$ given in \ref{o_convention}. Consider
the pair of o-planar representatives $(\Gg,\Gd(o_i))$ of $(\Gg,u)$
which respectively give  $(\Gt,o_i)$ after contracting their edges.

\begin{lem} \label{lem_reverse}
The o-planar tree $(\Gg,\Gd(o_1))$ is obtained by reversing the
o-planar structure $\Gd(o_2)_v$ of $(\Gg,\Gd(o_1))$ at vertex $v$
where $|F_\Gg(v)\bigcap \fF| \leq 1$.
\end{lem}

\begin{figure}[htb]
\centerfig{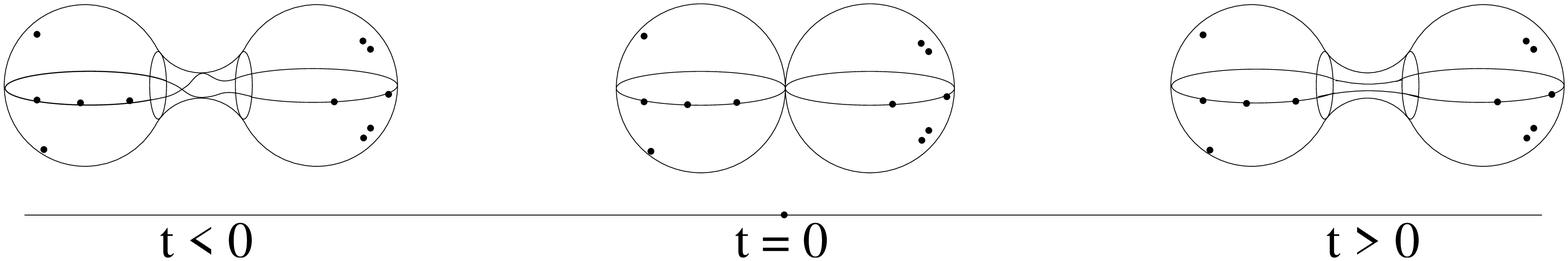,height=.8in} %
\caption{Two possible deformation of a real double point.}%
\label{fig_deform}
\end{figure}

\begin{proof}
Obviously, $(\Gg,\Gd(o_1))$ can be obtained from $(\Gg,\Gd(o_2))$ by
reversing the o-planar structures at one or both its vertices
$v_e,v^e$. If we reverse the o-planar structure of $(\Gg,\Gd(o_2))$
at the vertex $v$ such that $|F_\Gg(v)\bigcap \fF|> 1$, or at both
of its vertices $v_e$ and $v^e$, then the resulting o-planar tree
will not be an element of $\otree$ after contracting its edge:
reversing the o-planar structure at the vertex $v$ with
$|F_\Gg(v)\bigcap \fF|> 1$, or at both of the vertices reverses
cyclic order of the elements $\{2k+1,n-1,n\}$ when $l \geq 3$, and
moves $k$ from $\Conj^+$ to $\Conj^-$ when $l<3$ (see Figure
\ref{fig_deform}).
\end{proof}

For a pair of two-vertex o-planar trees $(\Gg,\Gd(o_i))$ as above,
we calculate the differences of parities as follows.

\begin{lem} \label{lem_parity}
Let $(\Gg,\Gd(o_i)), i=1,2$ be a pair of o-planar trees as above.
Let $V_\Gg=V_\Gg^\R = \{v_e,v^e\}$, and let o-planar structures at
the vertices $v_e$ and $v^e$ be
\begin{eqnarray*}
\begin{array}{lll}
\Gd(o_1)_{v_e} &=& \left\lbrace
\begin{array}{l}
\R\GS_{v_e} \ne \emp; F^\pm_{\Gg}(v_e); F_{\Gg}^\R (v_e) = \{
\{i_1\} < \cdots
< \{i_{q}\}< \{f_e\} < \\
<\{i_{q+r+1} \} < \cdots  <\{i_{l-1}\} <\{n\}   \}
\end{array}
\right\rbrace  \\
\Gd(o_2)_{v^e} &=& \{\R\GS_{v^e} \ne \emp; F^\pm_{\Gg}(v^e) ;
F_{\Gg}^\R (v^e) =\{ \{i_{q+1}\} < \cdots <\{i_{q+r}\}<\{f^e\} \}\}.
\end{array}
\end{eqnarray*}
Let $v$ be the vertex such that $|F_\Gg(v) \bigcap \fF| \leq 1$.
Then,
\begin{equation*}
|o_1|-|o_2|  = \left\lbrace
\begin{array}{ll}
 |F^+_{\Gg}(v^e)|+ \frac{r(r-1)}{2}& \mathrm{if}\ v = v^e \\
 |F^+_{\Gg}(v_e)|+ qr+rs+qs +\frac{s(s-1)}{2} + \frac{q(q-1)}{2}&
\mathrm{if}\ v = v_e \ \mathrm{and} \ |F_\Gg^\R(v_e)| > 3, \\
 |F^+_{\Gg}(v_e)|+ |F^\R_{\Gg}(v^e)|-1 &
\mathrm{if}\ v = v_e \ \mathrm{and} \ |F_\Gg^\R(v_e)| = 3, \\
 |F^+_{\Gg}(v_e)|  &
\mathrm{if}\ v = v_e \ \mathrm{and} \ |F_\Gg^\R(v_e)| = 2. \\
\end{array}
\right.
\end{equation*}
Here, $r = |F^\R_{\Gg}(v^e)|-1$ and $s = |F^\R_{\Gg}(v_e)|-q-2$.
\end{lem}

In Section \ref{sec_d_ori}, we have introduced differential forms
$\GO_{(\Gg_v,\Gd(o_i)_v)}$ for each $v \in V_\Gg$. When we reverse
the o-planar structure at the vertex $v$, the differential forms
$\GO_{(\Gg_v,\Gd(o_2)_v)}$, $\GO_{(\Gg_v,\Gd(o_1)_v)}$ become
related as follows.

\begin{lem} \label{lem_eksibir}
Let $(\Gg,\Gd(o_i)), i=1,2$ be two-vertex o-planar trees as above.
Then,
\begin{eqnarray*}
\GO_{(\Gg_v,\Gd(o_1)_v)} = (-1)^{\mu(v)} \ \GO_{(\Gg_v,\Gd(o_2)_v)},
\end{eqnarray*}
where
\begin{equation*}
\mu(v) = |F_{\Gg}^+(v)| +
\frac{(|F_{\Gg}^\R(v)|-2)(|F_{\Gg}^\R(v)|-3)}{2}.
\end{equation*}
\end{lem}

Lemmata \ref{lem_parity} and \ref{lem_eksibir} follow from
straightforward calculations.

\subsection{The first Stiefel-Whitney class}  \label{sec_sw_class}
This section is devoted to the proof of the following theorem.

\begin{thm}
\label{prop_sw1+}%

(i) For $\fix(\Gs) \ne \emp$, the Poincare dual of the  first Stiefel-Whitney class
of $\rmod{2k,l}$ is
\begin{eqnarray*}
[w_1] \ =\ \sum_{(\Gg,u)} \ [\csqc{\Gg,u}] \ =\ \sum_{\Gg}     \
[\real{\overline{D}_\Gg}] \ \mod 2,
\end{eqnarray*}
where the both sums are taken over all two-vertex trees such that
\begin{itemize}

\item $|F_\Gg(v^e) \bigcap \fF| \leq 1$ and  $|v^e| =0 \mod 2$, or

\item $|F_\Gg(v_e) \bigcap \fF| \leq 1$, $|F^\R_\Gg(v_e)|\ne 3$ and
$|v_e| (|v^e|-1) =0\mod 2$, or

\item $|F_\Gg(v_e) \bigcap \fF| \leq 1$, $|F^\R_\Gg(v_e)|  = 3$ and
$|F^\R_\Gg(v^e)|= 1$,
\end{itemize}
and, in the first sum, in addition over all $u$-planar structures on
$\Gg$.

(ii) For $\fix(\Gs) = \emp$, the Poincare dual of the  first Stiefel-Whitney
class of $\rmod{2k,0}$ vanishes.
\end{thm}

\begin{proof}
Fix an orientation for each top-dimensional stratum as in
\ref{o_convention}. The orientation $(-1)^{|o|} \ori{\Gt,o}$ of a
top-dimensional stratum $\csq{\Gt,o}$ induces some orientation of
each codimension one stratum $\csq{\Gg,\Gd(o)}$ (and
$\csq{\Gg,\hat{o}}$) contained in the boundary of $\csqc{\Gt,o}$.
The induced orientations $(-1)^{|o|} \ori{\Gg,\Gd(o)}$ and
$(-1)^{|o|}  \ori{\Gg,\hat{o}}$ are determined in Lemmata
\ref{lem_detay0}, \ref{lem_detay1} and  \ref{lem_detay2}, and they
give (relative) fundamental cycles $[\csqc{\Gg,\Gd(o)}]$ and
$[\csqc{\Gg,\hat{o}}]$ of the codimension one strata
$\csqc{\Gg,\Gd(o)}$ and $\csqc{\Gg,\hat{o}}$ respectively.

The Poincare dual of the first Stiefel-Whitney class of
$\rmod{2k,l}$ is given by
\begin{equation} \label{eqn_sw1}
[w_1] = \left\lbrace
\begin{array}{ll}
\frac{1}{2} \sum_{(\Gt,u) }  \left( \sum_{(\Gg,\Gd(o))}
[\csqc{\Gg,\Gd(o)}] \right) \ \mod 2, &
\mathrm{when}\ l>0, \\
\frac{1}{2} \sum_{(\Gt,u) }  \left( \sum_{(\Gg,\hat{o})}
[\csqc{\Gg,\hat{o}}] \right) \ \mod 2, & \mathrm{when}\ l=0,
\end{array}\right.
\end{equation}
where the external summation runs over all one-vertex u-planar trees
$(\Gt,u)$ and the internal one over all codimension one strata of
$\csqc{\Gt,o}$ for the one-vertex o-planar tree $(\Gt,o)$ which
represents $(\Gt,u)$ in accordance with \ref{o_convention}. Indeed,
the sum (\ref{eqn_sw1}) detects where the orientation on
$\rmodo{2k,l}$ can not be extended to $\rmod{2k,l}$.

We prove the theorem by evaluating  (\ref{eqn_sw1}).

\paragraph{\bf Case $\fix(\Gs) \ne \emp$.}
Let  $\csq{\Gt,o_i},i=1,2$ be a pair of adjacent top-dimensional
strata, and $\csq{\Gg,\Gd(o_i)} \subset \csqc{\Gt,o_i}$ be their
common codimension one boundary stratum. We calculate $\left[
\csqc{\Gg,\Gd(o_1)}\right]  + \left[ \csqc{\Gg,\Gd(o_2)}\right]$ as
follows. According to \ref{o_convention}, the strata $\csq{\Gt,o_i}$
are oriented by $(-1)^{|o_i|} \ori{\Gt,o_i}$, and these orientations
induce the orientations $(-1)^{|o_i|} \ori{\Gg,\Gd(o_i)}$ on
$\csq{\Gg,\Gd(o_i)}$. The induced orientations $(-1)^{|o_i|}
\ori{\Gg,\Gd(o_i)}$ are given by 
\begin{eqnarray*}
(-1)^{|o_i| + \GH_i}
[\GO_{(\Gg_{v^e},\Gd(o_i)_{v^e})} \wedge
\GO_{(\Gg_{v_e},\Gd(o_i)_{v_e})}]
\end{eqnarray*}
in Lemmata \ref{lem_detay0} and
\ref{lem_detay1} according to the convention introduced in Section
\ref{sec_d_ori}. We  denote by $v$ be the vertex such that
$|F_\Gg(v) \bigcap \fF| \leq 1$ as in Section \ref{sec_ad}, and
compare the induced orientations by calculating
\begin{equation*}
\Pi(o_1,o_2) =(|o_1|+\GH_1) -(|o_2| + \GH_2) - \mu(v)
\end{equation*}
for each of the following three subcases.

First, assume that $|F_\Gg(v^e) \bigcap \fF| \leq 1$. In this
subcase, the o-planar structure is reversed at the vertex $v=v^e$.
Therefore, $\GH_1 = \GH_2$ according to  Lemma \ref{lem_detay0}.
Finally, by applying Lemmata \ref{lem_parity} and \ref{lem_eksibir}
and using relation $r=|F_{\Gg}^\R(v^e)|-1$ we obtain
\begin{eqnarray*}
\Pi(o_1,o_2) =|o_1|-|o_2| - \mu(v^e) &=&  \frac{r(r-1)}{2}-
\frac{(|F_{\Gg}^\R(v^e)|-2)(|F_{\Gg}^\R(v^e)|-3)}{2} \\
&=& |F_{\Gg}^\R(v^e)|-2 \\
&=& |v^e| \mod 2.
\end{eqnarray*}
The latter equality follows from the fact that $|F_{\Gg}^\R(v)|= |v|
\mod 2$.

Second, assume that $|F_\Gg(v_e) \bigcap \fF| \leq 1$ and
$|F_\Gg^\R(v_e)| \ne 3$. In this subcase, the o-planar structure is
reversed at the vertex $v=v_e$. Since $|F_\Gg^\R(v_e)| \ne 3$, once
more $\GH_1 = \GH_2$ according to the Lemma \ref{lem_detay0}.
Finally, by applying Lemmata \ref{lem_parity} and \ref{lem_eksibir}
and using relation $|F_{\Gg}^\R(v_e)|=q+s+2$, we obtain
\begin{eqnarray*}
\Pi(o_1,o_2) &=&  qr+rs+qs +\frac{s(s-1)}{2} + \frac{q(q-1)}{2}
-\frac{(q+s)(q+s-1)}{2} \\
&=& r (q+s), \\
&=& (|F_{\Gg}^\R(v^e)|-1)(|F_{\Gg}^\R(v_e)|-2), \\
&=& |v_e|(|v^e|-1) \mod 2
\end{eqnarray*}
when $|F_\Gg^\R(v_e)| > 3$, and
\begin{eqnarray*}
\Pi(o_1,o_2) &=&  2 |F^+_\Gg(v_e)| = 0 \mod 2 \\
             &=& |v_e|(|v^e|-1) \mod 2
\end{eqnarray*}
when $|F_\Gg^\R(v_e)| =2$.

Third, we consider $|F_\Gg(v_e) \bigcap \fF| \leq 1$ and
$|F_\Gg^\R(v_e)| = 3$ case. In this subcase, the o-planar structure
is reversed at the vertex $v_e$. Hence, $\GH_1 = \GH_2$ whenever
$|F_\Gg^\R(v^e)| = 1,2$, and $\GH_1 - \GH_2$ is $\pm(l+1)=\pm
(|F_{\Gg}^\R(v^e)|+2)$ whenever $|F_\Gg^\R(v^e)| \geq 3$. Finally,
by applying Lemmata \ref{lem_parity} and \ref{lem_eksibir}, we
obtain
\begin{eqnarray*}
\Pi(o_1,o_2) = \left\lbrace
\begin{array}{lll}
|F_{\Gg}^\R(v^e)|-1 \pm (|F_{\Gg}^\R (v^e)|+2)
& = 1 \mod 2, & \mathrm{when} \ |F_\Gg^\R(v^e)| \geq 3, \\
|F_{\Gg}^\R(v^e)|-1
& = 1 \mod 2, & \mathrm{when} \ |F_\Gg^\R(v^e)| = 2, \\
|F_{\Gg}^\R(v^e)|-1 & = 0 \mod 2, & \mathrm{when} \ |F_\Gg^\R(v^e)|
= 1,
\end{array}
\right.
\end{eqnarray*}

The induced orientations $(-1)^{|o_i|} \ori{\Gg,\Gd(o_i)}$ are the
same if and only if $\Pi(o_1,o_2) = 0 \mod 2$. Hence, we have
\begin{eqnarray*}
\left[ \csqc{\Gg,\Gd(o_1)}\right]  + \left[
\csqc{\Gg,\Gd(o_2)}\right] &=&  \frac{1 + (-1)^{\Pi(o_1,o_2)}}{2}
\left[ \csqc{\Gg,\Gd(o_1)}\right].
\end{eqnarray*}

The sum $([\csqc{\Gg,\Gd(o_1)}] + [\csqc{\Gg,\Gd(o_2)}])/2$ gives us
the fundamental cycle $[\csqc{\Gg,\Gd(o_1)}]$ when
$(-1)^{|o_1|}\ori{\Gg,\Gd(o_1)}= (-1)^{|o_2|}\ori{\Gg,\Gd(o_2)}$,
and it turns to zero otherwise. Finally, as is follows from the
above case-by-case calculations of $\Pi(o_1,o_2)$, the fundamental
class of a codimension one strata $\csqc{\Gg,\Gd(o_1)}$ is involved
in $[w_1]$ if and only if one of the three conditions given in
Theorem are verified. It gives the first expression for $[w_1]$
given in Theorem. Since in this first expression the sum is taken
over all u-planar structures on $\Gg$, it can be shorten to the sum
of the fundamental classes of $\real{D_\Gg}$.

\paragraph{\bf Case  $\fix(\Gs) = \emp$.}
Let $\csq{\Gt,u_i},i=1,2,$ be a pair of adjacent top-dimensional
strata and $\csq{\Gg,u}$ be their common codimension one boundary
stratum. Let $(\Gt,o_i)$ be the o-planar representatives of
$(\Gt,u_i)$ given in \ref{o_convention}. Here,  we have to consider
two subcases: (i) $\csq{\Gg,u}$ is a stratum of  real curves with
two real components (i.e., $|V_\Gg|=|V_\Gg^\R| =2$), and (ii)
$\csq{\Gg,u}$ is a stratum of real curves with two complex
conjugated components (i.e., $|V_\Gg|=2$ and $|V_\Gg^\R| =0$).

(i) Consider the pair of o-planar representatives $(\Gg,\Gd(o_i))$
of $(\Gg,u)$ which respectively give  $(\Gt,o_i)$ after contracting
the edges and compare their o-planar structure. Since the both tails
$n$ and $\Gs(n)$ are in $F_\Gg(v_e)$, the o-planar structure is
reversed at the vertex $v^e$. Therefore, $\GH_1 = \GH_2$ according
to the Lemma \ref{lem_detay0}. Finally, by applying Lemmata
\ref{lem_parity} and \ref{lem_eksibir}, we obtain
\begin{eqnarray*}
\Pi(o_1,o_2) = 2 |F_\Gg^+(v^e)| +1 = 1 \mod 2.
\end{eqnarray*}
In other words, $\left[ \csqc{\Gg,\Gd(o_1)}\right]  + \left[
\csqc{\Gg,\Gd(o_2)}\right]  = 0$ for this case.

(ii) Let $\csq{\Gt,o_2}$  be a stratum of real curves with empty
real part, and let $(\Gg,\hat{o})$ be an o-planar representative of
$(\Gg,u)$.

The orientations of $\csq{\Gg,u}$ induced by the  orientations
$(-1)^{|o_1|}\ori{\Gt,o_1}$ and $\ori{\Gt,o_2}$ of $\csq{\Gt,o_1}$
and $\csq{\Gt,o_2}$ are given in Lemmata \ref{lem_detay1} and
\ref{lem_detay2}. Namely, they are respectively given by the
following differential forms
\begin{eqnarray*}
(-1)^{|o_1|} \bigwedge_{\Ga_* \in F_{\Gg}^+ \smin (F_{\Gg}^+ \bigcap
\{k-1,k,2k-1,2k\}) }
dz_{\Ga_*} \wedge d\bar{z}_{\Ga_*}, \\
(-1)^{\aleph} \bigwedge_{\Ga_* \in F_{\Gg}^+ \smin (F_{\Gg}^+
\bigcap \{k-1,k,2k-1,2k\}) }
dz_{\Ga_*} \wedge d\bar{z}_{\Ga_*}, \\
\end{eqnarray*}
where $|o_1|=|\{1,\cdots,k-1\} \bigcap F_\Gg^-|$ and $\aleph =
|\{1,\cdots,k-1\} \bigcap F_\Gg^-|+1$. Therefore, the orientations
induced from different sides are opposite and the sum
$(-1)^{\aleph-1} \left[ \csqc{\Gg,u}\right] + (-1)^{\aleph} \left[
\csqc{\Gg,u}\right]$ vanishes for all such $(\Gg,\hat{o})$.
\end{proof}

\subsubsection{Example}
\begin{figure}[htb]
\centerfig{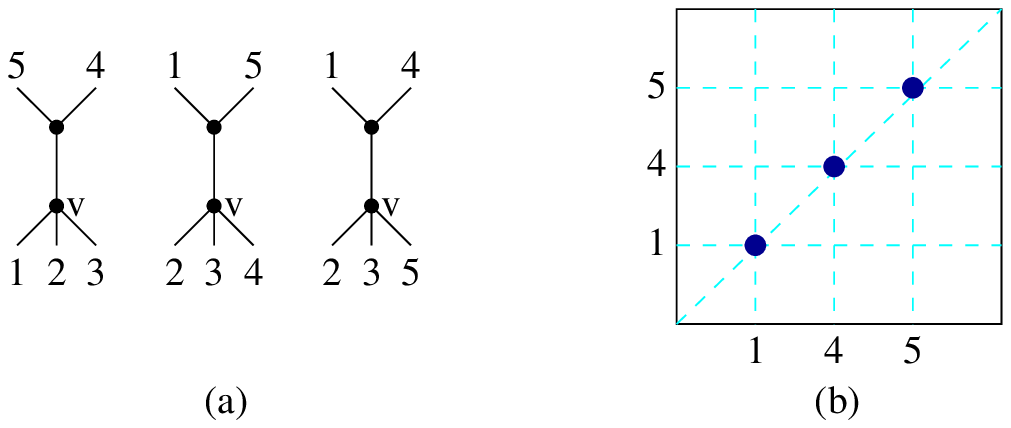,height=4cm} %
\caption{(a) The $5$-trees of Stiefel-Whitney class of $\rmod{0,5}$
due Theorem \ref{prop_sw1+},
(b) The blown-up locus in $\rmod{0,5}$} %
\label{fig_disorienttree}
\end{figure}

Due to Theorem \ref{prop_sw1+}, the Poincare dual of the first
Stiefel-Whitney class $[w_1]$ of $\rmod{0,5}$ can be  represented by
$\sum_\Gg [\real{\overline{D}_{\Gg}}] = \sum_{(\Gg,u)}
[\csqc{\Gg,u}]$ where $\Gg$ are $5$-trees with a vertex $v$
satisfying $|v|=4$ and $|F_\Gg(v) \bigcap \{1,4,5\}|=1$. These
$5$-trees are given in Figure \ref{fig_disorienttree}a, and the
union corresponding strata $\bigcup_\Gt \real{\overline{D}_{\Gg}}$
is given the three exceptional divisors obtained by blowing up the
three highlighted points in Figure \ref{fig_disorienttree}b.

\section{The orientation double covering of  $\rmod{2k,l}$}
\label{sec_coverofmoduli}

In Section \ref{sec_sw_class}, the first Stiefel-Whitney class of
$\rmod{2k,l}$ is determined in terms of its strata. We have also
proved that the moduli space $\rmod{2k,l}$ is orientable when $n=4$
or $l=0$. In this section, we give a combinatorial construction of
orientation double covering for the rest of the cases i.e., $n>4$ and
$l>0$. By observing the non-triviality of the orientation double
cover in these cases, we prove that $\rmod{2k,l}$ is not orientable.

\subsection{Construction of orientation double covering}
\label{sec_glue}

In Section \ref{sec_conf_ir}, we have shown that the map
$\cspq{2k,l} \to \rmodo{2k,l}$, which is identifying the reverse
o-planar structures, is a trivial double covering. The disjoint
union of closed strata $\csqc{2k,l} = \bigsqcup_{(\Gt,o)}
\csqc{\Gt,o}$, where $|V_\Gt|=1$ and $(\Gt,o)$ runs over all
possible o-planar structures on $\Gt$, is a natural compactification
of $\cspq{2k,l}$.

To obtain the orientation double covering of $\rmod{2k,l}$ we need
to get rid of the codimension one strata by pairwise gluing them. We
use the following simple recipe: for each pair $(\Gt,o_i), =1,2,$ of
one-vertex o-planar trees obtained by contracting the edge in a pair
$(\Gt,\Gd(o_i)), i=1,2,$ of two-vertex o-planar trees with the same
underlying tree such that $V_\Gg=V_\Gg^\R = \{v_e,v^e\}$, $v_e=
\dd_\Gg(n)$, we glue $\csqc{\Gt,o_i}$ along
$\csqc{\Gg,\Gd(o_i)},i=1,2$, if
\begin{itemize}

\item[$\cA$.] $(\Gg,\Gd(o_1))$ produces $(\Gg,\Gd(o_2))$ by reversing the
o-planar structure at the vertex $v^e$,  $|F_\Gg(v^e) \bigcap \fF|
\leq 1$, and $|v^e|=1 \mod 2$,

\item[$\cB$.] $(\Gg,\Gd(o_1))$ produces $(\Gg,\Gd(o_2))$ by reversing the
o-planar structure at the vertex $v_e$,  $|F_\Gg(v^e) \bigcap \fF|
\leq 1$, and $|v^e|=0 \mod 2$,

\item[$\cC$.] $(\Gg,\Gd(o_1))$ produces $(\Gg,\Gd(o_2))$ by reversing the
o-planar structure at the vertex $v_e$, $|F_\Gg(v_e) \bigcap \fF|
\leq 1$,
 $|F_\Gg^\R(v_e)| \ne 3$ and $|v_e|(|v^e-1|)=1 \mod 2$,

\item[$\cD$.] $(\Gg,\Gd(o_1))$ produces $(\Gg,\Gd(o_2))$ by reversing the
o-planar structure at the vertex $v^e$, $|F_\Gg(v_e) \bigcap \fF|
\leq 1$,
 $|F_\Gg^\R(v_e)| \ne 3$ and $|v_e|(|v^e-1|)=0 \mod 2$,

\item[$\cE$.] $(\Gg,\Gd(o_1))$ produces $(\Gg,\Gd(o_2))$ by reversing the
o-planar structure at the vertex $v_e$, $|F_\Gg(v_e) \bigcap \fF|
\leq 1$,
 $|F_\Gg^\R(v_e)| = 3$ and $|F_\Gg^\R(v^e)| \ne 1$,

\item[$\cF$.] $(\Gg,\Gd(o_1))$ produces $(\Gg,\Gd(o_2))$ by reversing the
o-planar structure at the vertex $v^e$, $|F_\Gg(v_e) \bigcap \fF|
\leq 1$,
 $|F_\Gg^\R(v_e)| = 3$ and $|F_\Gg^\R(v^e)|  = 1$.

\end{itemize}

We denote by $\cover{2k,l}$ the resulting factor space.

\begin{thm}
\label{thm_cover} %
$\cover{2k,l}$ is the orientation double cover of $\rmod{2k,l}$.
\end{thm}

\begin{proof}
Let $\ove{M}$ be the orientation double covering of $\rmod{2k,l}$.
The points of $\ove{M}$ can be considered as points in $\rmod{2k,l}$
with local orientation. On the other hand, by using opposite
o-planar structures of a one-vertex $\Gt$  we can determine
orientations $(-1)^{|o|} \ori{\Gt,o}$ and $(-1)^{|\bar{o}|+l-1}
\ori{\Gt,\bar{o}}$ on $\csq{\Gt,o}$ and $\csq{\Gt,\bar{o}}$ where
$(\Gt,o) \in \otree$. These orientations are opposite with respect
to the identification of $\csq{\Gt,o}$ and $\csq{\Gt,\bar{o}}$ by
the canonical diffeomorphism  $-\I$ introduced in Subsection
\ref{sec_conf_ir}. Hence, there is a natural continuous embedding
$\cspq{2k,l} =\bigsqcup_{(\Gt,o)\in  \otree } (\csq{\Gt,o}\sqcup
\csq{\Gt,\bar{o}})\to \ove{M}$. It extends to a surjective
continuous map $\csqc{2k,l} = \bigsqcup_{(\Gt,o)\in  \otree }
(\csqc{\Gt,o}\sqcup \csqc{\Gt,\bar{o}}) \to \ove{M}.$ Since
$\csqc{2k,l}$ is compact and $\ove{M}$ is Hausdorff, the orientation
double covering $\ove{M}$ is a quotient space
$\csqc{2k,l}/\textsl{R}$ of $\csqc{2k,l}$ under the equivalence
relation $\textsl{R}$ defined by the map  $\csqc{2k,l} \to \ove{M}$.

This equivalence relation is uniquely determined by its restriction
to the codimension one faces of $\csqc{2k,l}$, which cover the
codimension one strata of $\rmod{2k,l}$ under the composed map
$\csqc{2k,l} \to \ove{M}\to \rmod{2k,l}$. On the other hand, the
equivalence relation on the codimension one faces is determined by
the first Stiefel-Whitney class: A partial section of the induced
map  $\csqc{2k,l}/\textsl{R} \to \rmod{2k,l}$ given by distinguished
strata $\bigsqcup_{(\Gt,o)\in  \otree }  \csq{\Gt,o}$. Over a
neighborhood of a codimension one stratum of $\rmod{2k,l}$, a
partial section extends to a section if this strata is not involved
in the expression for the first Stiefel-Whitney class given in
Theorem \ref{prop_sw1+}, and it should not extend, otherwise. Notice
that the faces $\csqc{\Gt,\Gd(o_i)}$ considered in relations $\cA$,
$\cC$ and $\cE$ are mapped onto the strata $\csqc{\Gt,\Gd(u)}$ which
do not contribute to the expression $[w_1]$ given in Theorem
\ref{prop_sw1+}, and the faces $\csqc{\Gt,\Gd(o_i)}$  in relations
$\cB$, $\cD$ and $\cF$ are mapped onto the strata
$\csqc{\Gt,\Gd(u)}$ which contribute to the expression $[w_1]$.
There are four different faces  $\csqc{\Gt,\Gd(o)_i}, i=1,\cdots,4$ over
each codimension one stratum $\csqc{\Gt,\Gd(u)}$.  Lemma
\ref{lem_reverse} determines the pairs $\csqc{\Gt,\Gd(o)_i}$,
$\csqc{\Gt,\Gd(o)_j}$ to be glued to  each other.
\end{proof}

\begin{cor}
The moduli space $\rmod{2k,l}$ is not orientable when $2k+l>4$ and
$l>0$.
\end{cor}

\begin{proof}
Let $l \geq 3$, and  $(\Gt,o)$ be an o-planar structure with
$\{2k+1\}<\{n-1\} <\{n\}$. It is clear that, we can produce any
o-planar structure on $\Gt$ with $\{2k+1\}<\{n-1\} <\{n\}$ by
applying following operations consecutively:
\begin{itemize}
\item  interchanging the order  of two consecutive
tails $\{i,i+1\}$ for $|\{i,i+1\} \bigcap \fF| \leq 1$  and $n
\not\in \{i,i+1\}$,

\item swapping $j \in \Conj^+$ with $\bar{j}\in
\Conj^-$ for $j \ne k,{2k}$.
\end{itemize}
The one-vertex o-planar trees with $\{n-1\}<\{2k+1\} <\{n\}$  can be
produced from the o-planar tree $(\Gt,\bar{o})$ via same procedure.

Let $l=1,2$. Similarly, if  we start with o-planar tree $(\Gt,o)$
with $k \in \Conj^+$ ($k \in \Conj^-$), we can produce any o-planar
structure on $\Gt$ with $k \in \Conj^+$ ($k \in \Conj^-$) by
swapping $j \in \Conj^+$ with $\bar{j}\in \Conj^-$ for $j \ne
k,{2k}$.

Note that, these operations correspond to passing from one
top-dimensional stratum to another in $\cover{2k,l}$  through
certain faces. These faces correspond to the one-edge o-planar trees
$(\Gg,\Gd(o)_i)$ with $F_\Gg(v^e)=\{i,{i+1},f^e\}$ (resp.
$F_\Gg(v^e)=\{j,\bar{j},f^e\}$) which
 are faces  glued according to the relations
of type $\cA$. Hence, any two top-dimensional strata in
$\cover{2k,l}$ with same cyclic ordering of $\fF$ (resp. with $k$ is
in same set $\Conj^\pm$) can be connected through a path passing
through these codimension faces $\csqc{\Gg,\Gd(o)_i}$. The quotient
space $\csqc{2k,l}/\cA$ has two connected components since there are
two possible cyclic orderings of $\fF$   when $l\geq 3$ (resp. two
possibilities for $l=1,2$ case: $k \in \Conj^+$ and $k \in
\Conj^-$).

The set of relations of type $\cB$ is not empty when $2k+l>4$ and
$l>0$. Moreover, the relation of type $\cB$ reverses the cyclic
ordering on $\fF$ (resp. moves $k$ from $\Conj^\pm$ to $\Conj^\mp$).
Hence, the faces glued according to the relations of type $\cB$
connect two different components of $\csqc{2k,l}/\cA$. Therefore,
the orientation double cover $\rmod{2k,l}$ is nontrivial when
$2k+l>4$ and $l>0$ which simply means that the moduli space
$\rmod{2k,l}$ is not orientable in this case.
\end{proof}

\subsubsection{Examples}
\begin{figure}[htb]
\centerfig{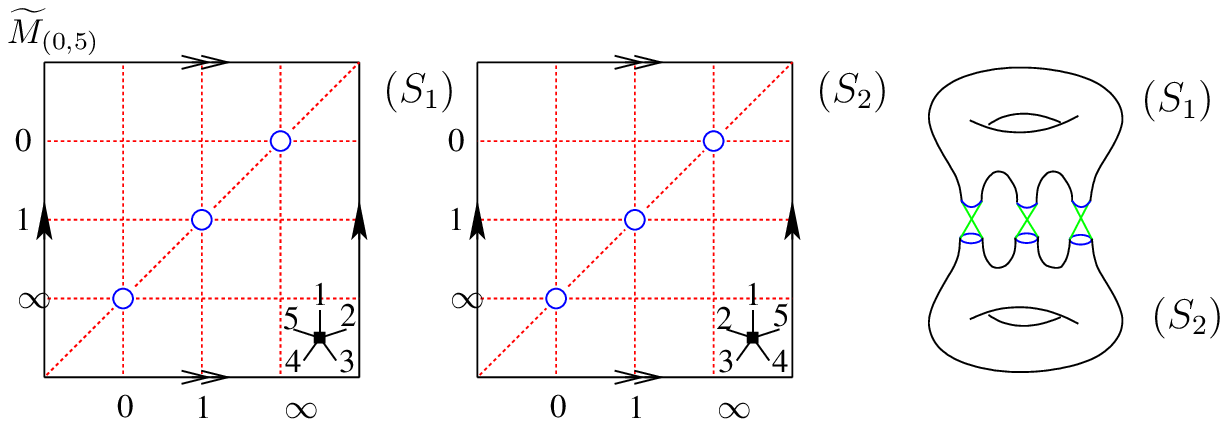,height=5cm} %
\caption{Stratification of $\cover{0,5}$} %
\label{fig_cover3}
\end{figure}

In Example \ref{exa_realmoduli}, we obtained that $\rmod{0,5},
\rmod{2,3}$, and $\rmod{4,1}$ are respectively, a torus with three
points blown up, a sphere with three points blown up, and a sphere
with one point blown up. The coverings $\cover{0,5},\cover{2,3}$ and
$\cover{4,1}$ are obtained by taking the two copies of the
corresponding moduli space of real curves and replacing the blown up
loci by annuli. Therefore, $\cover{0,5},\cover{2,3}$ and
$\cover{4,1}$ are surfaces of genus 4, genus 2 and genus 0,
respectively (see Figure \ref{fig_cover3} which illustrates the case
$(k,l)=(0,5)$).

\subsection{Combinatorial types  of strata of $\cover{2k,l}$}
\label{sec_r_equi}%
While constructing $\cover{2k,l}$, the closure of the each
codimension one strata are glued in a consistent way. This
identification of codimension strata gives an equivalence relation
among the o-planar trees when $l \not=0$.

We define the notion of {\it R-equivalence} on the set of such
o-planar trees by treating  different cases separately. Let
$(\Gg_1,o_1), (\Gg_2,o_2)$ be o-planar trees.
\begin{enumerate}

\item If $|V_{\Gg_i}^\R| = 1$, then we say that they are
R-equivalent whenever $\Gg_1,\Gg_2$ are isomorphic (i.e., $\Gg_1
\approx \Gg_2$) and  the o-planar structures are the same.

\item If $\Gg_i$ have an edge corresponding to real node
(i.e. $E_{\Gg_i}^\R =\{e\}$ and $V_{\Gg_i}^ \R= \dd_{\Gg}(e) =
\{v^e,v_e\}$), we first obtain $(\Gg_i (e),o_i(e))$ by contracting
conjugate pairs of edges until there will be none. We say that
$(\Gg_1,o_1)$ and $(\Gg_2,o_2)$ are R-equivalent whenever $\Gg_1
\approx \Gg_2$ and
\begin{itemize}
\item $(\Gg_1(e),o_1(e))$ produces $(\Gg_1(e),o_2(e))$ by reversing the
o-planar structure at the vertex $v^e$,  $|F_\Gg(v^e) \bigcap \fF|
\leq 1$, and $|v^e|=1 \mod 2$,

\item $(\Gg_1(e),o_1(e))$ produces $(\Gg_1(e),o_2(e))$ by reversing the
o-planar structure at the vertex $v_e$,  $|F_\Gg(v^e) \bigcap \fF|
\leq 1$, and $|v^e|=0 \mod 2$,

\item $(\Gg_1(e),o_1(e))$ produces $(\Gg_1(e),o_2(e))$ by reversing the
o-planar structure at the vertex $v_e$, $|F_\Gg(v_e) \bigcap \fF|
\leq 1$, $|F_\Gg^\R(v_e)| \ne 3$ and $|v_e|(|v^e-1|)=1 \mod 2$,

\item $(\Gg_1(e),o_1(e))$ produces $(\Gg_1(e),o_2(e))$ by reversing the
o-planar structure at the vertex $v^e$, $|F_\Gg(v_e) \bigcap \fF|
\leq 1$, $|F_\Gg^\R(v_e)| \ne 3$ and $|v_e|(|v^e-1|)=0 \mod 2$,

\item $(\Gg_1(e),o_1(e))$ produces $(\Gg_1(e),o_2(e))$ by reversing the
o-planar structure at the vertex $v_e$, $|F_\Gg(v_e) \bigcap \fF|
\leq 1$, $|F_\Gg^\R (v_e)| = 3$ and $|F_\Gg^\R (v^e)| \ne 1$,

\item $(\Gg_1(e),o_1(e))$ produces $(\Gg_1(e),o_2(e))$ by reversing the
o-planar structure at the vertex $v^e$, $|F_\Gg(v_e) \bigcap \fF|
\leq 1$, $|F_\Gg^\R (v_e)| = 3$ and $|F_\Gg^\R (v^e)| = 1$,
\end{itemize}

\item Otherwise, if $\Gg_i$ have more than one invariant edge
(i.e. $|E_{\Gg_i}|^\R>1$), we say that $(\Gg_1,o_1)$, $(\Gg_2,o_2)$
are R-equivalent whenever $\Gg_1 \approx \Gg_2$ and there exists an
edge $e \in E_{\Gg_i}^\R$ such that the o-planar trees
$(\Gg_{i}(e),o_i(e))$, which are obtained by contracting all edges
but $e$, are R-equivalent in the sense of the Case (2).
\end{enumerate}

We call the maximal set of pairwise R-equivalent o-planar trees by
{\it R-equivalence classes} of o-planar trees.

\begin{thm}
A stratification of the orientation double cover $\cover{2k,l}$ is
given by
\begin{eqnarray*}
\cover{2k,l} = \bigsqcup_{{\mathrm{R-equivalance \ classes }} \atop
{\mathrm{of \ o-planar} \ (\Gg,o)}} \csq{\Gg,o}.
\end{eqnarray*}
\end{thm}

\subsection{Some other double coverings of $\rmod{2k,l}$}

In \cite{ka}, Kapranov constructed a different double covering
$\widehat{\R M}_{(0,l)}$ of $\rmod{0,l}$ having no boundaries. He
has applied the following recipe to obtain the double covering: Let
$\overline{C}_{(0,l)}$ be the disjoint union of closed strata
$\bigsqcup_{(\Gt,o)}  \csqc{\Gt,o}$ as above. Let $(\Gg,\Gd(o_i)),
i=1,2$ be two-vertex o-planar trees representing the same u-planar
tree  $(\Gg,u)$, and let $(\Gt,o_i)$ be the one-vertex trees
obtained by contracting the edges of $(\Gg,\Gd(o_i))$. The strata
$\csqc{\Gg,\Gd(o_i)},i=1,2$ are glued if $(\Gg,\Gd(o_1))$ produces
$(\Gg,\Gd(o_2))$ by reversing the o-planar structure at vertex $v^e
\ne \dd_\Gg(n)$. We obtain first Stiefel-Whitney class of
$\widehat{\R M}_{(0,l)}$ by using the same arguments in Theorem
\ref{prop_sw1+}.

\begin{prop}
The Poincare dual of the  first Stiefel-Whitney class of
$\widehat{\R M}_{(0,l)}$ is
\begin{eqnarray*}
[\widehat{w_1}] = \frac{1}{2}
\sum_{(\Gt,o)} \sum_{(\Gg,\Gd(o)): |v^e| =0 \mod 2}
 \left[ \csqc{\Gg,\Gd(o)}\right] \mod 2.
\end{eqnarray*}
\end{prop}

It is well-known that these spaces are not orientable when $l \geq
5$.

\subsubsection{A double covering from open-closed string theory}

In \cite{fu,liu}, a different `orientation double covering' is
considered. It can be given as the disjoint union
$\bigsqcup_{(\Gt,o)} \csqc{\Gt,o}$ where $F_{\Gt}^+ =
\{1,\cdots,k\}$, and $F_{\Gt}^\R$ carries all possible oriented
cyclic ordering. It is a disjoint union of manifolds with corners.
The covering map   $\bigsqcup_{(\Gt,o)} \csqc{\Gt,o} \to
\rmod{2k,l}$ is two-to-one only over a subset of the open space
$\rmodo{2k,l}$. It only covers the subset $\bigsqcup_{(\Gt,u)}
\csqc{\Gt,u}$ of $\rmod{2k,l}$ where  u-planar trees $(\Gt,u)$ have
the partition $\{\{1,\cdots,k\},\{k+1,\cdots,2k\}\}$ of $F_\Gt \smin
F_\Gt^\R$. Moreover, the covering map is not two-to-one over the
strata with codimension higher than zero.

{\bf Acknowledgements:} I take this opportunity to express my deep
gratitude to my supervisors V. Kharlamov and M. Polyak for their
assistance and guidance in this work. I also wish to thank to Yu. I.
Manin for his interest and suggestions which have been invaluable
for this work.

I am thankful to  K. Aker,  E. Ha, A. Mellit, D. Radnell and A. Wand
for their comments and suggestions.


Thanks are also due to Max-Planck-Institut f\"ur Mathematik, Israel
Institute of Technology and Institut de Recherhe Math\'{e}matique
Avanc\'{e}e de Strasbourg for their hospitality.

\end{document}